\documentclass[12pt,reqno]{amsart}

\usepackage[arrow,matrix,curve]{xy}

\usepackage[dvips]{graphicx}

\usepackage{amssymb, latexsym, amsmath, amscd, array,
}

\DeclareMathOperator{\adequal}{\,{}_{\ulcorner\!\urcorner}\,}

\theoremstyle{definition}

\numberwithin{equation}{section}

\newcommand\N {{\mathbb N}}
\newcommand\R {{\mathbb R}}

\newcommand\RRR{\mathbb{I\hskip-3pt R}}

\newcommand\Los{{\L}o{\'s}}

\newcommand\parisotes{{$\pi\alpha\rho\iota\sigma\acute{o}\tau\eta\varsigma$}}

\newcommand\megethos
{{$\mu\acute\varepsilon\gamma\varepsilon\vartheta{}o\varsigma$}}

\author[J.B.]{Jacques Bair}\address{J. Bair, HEC-ULG, University of
Liege, 4000 Belgium}\email{j.bair@ulg.ac.be}

\author[P.B.]{Piotr B\l{}aszczyk}\address{P. B\l{}aszczyk, Institute
of Mathematics, Pedagogical University of Cracow,
Poland}\email{pb@up.krakow.pl}

\author[R.E.]{Robert Ely}\address{R. Ely, Department of Mathematics,
University of Idaho, Moscow, ID 83844 US}\email{ely@uidaho.edu}

\author[V.H.]{Val\'erie Henry}\address{V.~Henry, Department of
Mathematics, University of Namur, 5000
Belgium}\email{vhen@math.fundp.ac.be}

\author[V.K.]{Vladimir Kanovei} \address{V. Kanovei, IPPI, Moscow, and
MIIT, Moscow, Russia} \email{kanovei@rambler.ru}

\author[K.K.]{Karin U. Katz}\address{K. Katz, Department of
Mathematics, Bar Ilan University, Ramat Gan 52900 Israel}

\author[M.K.]{Mikhail G. Katz}\address{M. Katz, Department of
Mathematics, Bar Ilan University, Ramat Gan 52900 Israel}
\email{katzmik@macs.biu.ac.il}

\author[S.K.]{Semen S. Kutateladze}\address{S. Kutateladze, Sobolev
Institute of Mathematics, Novosibirsk State University, Russia}
\email{sskut@math.nsc.ru}

\author[T.M.]{Thomas McGaffey}\address{T. McGaffey, Rice University,
US}\email{thomasmcgaffey@sbcglobal.net}

\author[D.Sc.]{David M. Schaps}\address{D. Schaps, Department of
Classical Studies, Bar Ilan University, Ramat Gan 52900 Israel}
\email{dschaps@mail.biu.ac.il}

\author[D.Sh.]{David Sherry}\address{D. Sherry, Department of
Philosophy, Northern Arizona University, Flagstaff, AZ 86011,
US}\email{david.sherry@nau.edu}

\author[S.S.]{Steven Shnider}\address{S.~Shnider, Department of
Mathematics, Bar Ilan University, Ramat Gan 52900
Israel}\email{shnider@macs.biu.ac.il}

\begin{document}

\title{Is mathematical history written by the victors?}

\keywords{Adequality; Archimedean axiom; infinitesimal; law of
continuity; mathematical rigor; standard part principle;
transcendental law of homogeneity; variable quantity}

\begin{abstract}
We examine prevailing philosophical and historical views about the
origin of infinitesimal mathematics in light of modern infinitesimal
theories, and show the works of Fermat, Leibniz, Euler, Cauchy and
other giants of infinitesimal mathematics in a new light.  We also
detail several procedures of the historical infinitesimal calculus
that were only clarified and formalized with the advent of modern
infinitesimals.  These procedures include Fermat's adequality;
Leibniz's law of continuity and the transcendental law of homogeneity;
Euler's principle of cancellation and infinite integers with the
associated infinite products; Cauchy's infinitesimal-based definition
of continuity and ``Dirac'' delta function.  Such procedures were
interpreted and formalized in Robinson's framework in terms of
concepts like microcontinuity (S-continuity), the standard part
principle, the transfer principle, and hyperfinite products.  We
evaluate the critiques of historical and modern infinitesimals by
their foes from Berkeley and Cantor to Bishop and Connes.  We analyze
the issue of the consistency, as distinct from the issue of the rigor,
of historical infinitesimals, and contrast the methodologies of
Leibniz and Nieuwentijt in this connection.
\end{abstract}

\maketitle

\noindent

\tableofcontents

\section{The ABC's of the history of infinitesimal mathematics}
\label{one}

The
%
%
ABCs of the history of infinitesimal mathematics are in need of
clarification.  To what extent does the famous dictum ``history is
always written by the victors" apply to the history of mathematics, as
well?  A convenient starting point is a remark made by Felix Klein in
his book \emph{Elementary mathematics from an advanced standpoint}
(Klein 1908 \cite[p.~214]{Kl08}).  Klein wrote that there are not one
but two separate tracks for the development of analysis:

\begin{enumerate}
\item[(A)] the Weierstrassian approach (in the context of an
\emph{Archimedean} continuum); and
\item[(B)] the approach with indivisibles and/or infinitesimals (in
the context of what we will refer to as a \emph{Bernoullian}
continuum).%
\footnote{Systems of quantities encompassing infinitesimal ones were
used by Leibniz, Bernoulli, Euler, and others.  Our choice of the term
is explained in entry~\ref{be}.  It includes modern non-Archimedean
systems.}
\end{enumerate}

Klein's sentiment was echoed by the philosopher G. Granger, in the
context of a discussion of Leibniz, in the following terms:
\begin{quote}
Aux yeux des d\'etracteurs de la nouvelle Analyse, l'insur\-montable
difficult\'e vient de ce que de telles pratiques font violence aux
r\`egles ordinaires de l'Alg\`ebre, tout en conduisant \`a des
r\'esultats, exprimables en termes finis, dont on ne saurait contester
l'exactitude. Nous savons aujourd'hui que deux voies devaient s'offrir
pour la solution du probl\`eme:

[A] Ou bien l'on \'elimine du langage math\'ematique le terme
d'infiniment petit, et l'on \'etablit, en termes finis, le sens \`a
donner \`a la notion intuitive de `valeur limite'\ldots

[B] Ou bien l'on accepte de maintenir, tout au long du Calcul, la
pr\'esence d'objets portant ouvertement la marque de l'infini, mais en
leur conf\'erant un statut propre qui les ins\`ere dans un syst\`eme
dont font aussi partie les grandeurs finies\ldots

C'est dans cette seconde voie que les vues philosophi\-ques de Leibniz
l'ont orient\'e (Granger 1981 \cite[pp.~27-28]{Gran}).%
\footnote{Similar views were expressed by M.~Parmentier in (Leibniz
1989 \cite[p.~36, note~92]{Le89}).}
\end{quote}

Thus we have two parallel tracks for conceptualizing infinitesimal
calculus, as shown in Figure~\ref{AB}.

\begin{figure}
\[
\xymatrix@C=95pt{{} \ar@{-}[rr] \ar@{-}@<-0.5pt>[rr]
\ar@{-}@<0.5pt>[rr] & {}
& \hbox{\quad
B-continuum} \\ {} \ar@{-}[rr] & {} & \hbox{\quad A-continuum} }
\]
\caption{\textsf{Parallel tracks: a thick continuum and a thin
continuum}}
\label{AB}
\end{figure}

At variance with Granger's appraisal, some of the literature on the
history of mathematics tends to assume that the A-approach is the
ineluctably ``true" one, while the infinitesimal B-approach was, at
best, a kind of evolutionary dead-end or, at worst, altogether
inconsistent.  To say that infinitesimals provoked passions would be
an understatement.  Parkhurst and Kingsland, writing in \emph{The
Monist}, proposed applying a \emph{saline solution} (if we may be
allowed a pun) to the problem of the infinitesimal:
\begin{quote}
[S]ince these two words [infinity and infinitesimal] have sown
nearly as much faulty logic in the fields of mathematics and
metaphysics as all other fields put together, they should be rooted
out of both the fields which they have contaminated.  And not only
should they be rooted out, lest more errors be propagated by them:
\emph{a due amount of salt} should be ploughed under the infected
territory, that the damage be mitigated as well as arrested (Parkhurst
and Kingsland 1925 \cite[pp.~633-634]{PK}) [emphasis added--the
authors].
\end{quote}
Writes P.~Vickers:
\begin{quote}
So entrenched is the understanding that the early calculus was
inconsistent that many authors don't provide a reference to support
the claim, and don't present the set of inconsistent propositions they
have in mind.  (Vickers 2013 \cite[section~6.1]{Vi13})
\end{quote}

Such an \emph{assumption} of inconsistency can influence one's
appreciation of historical mathematics, make a scholar myopic to
certain significant developments due to their automatic placement in
an ``evolutionary dead-end" track, and inhibit potential fruitful
applications in numerous fields ranging from physics to economics (see
Herzberg 2013 \cite{Her}).  One example is the visionary work of
Enriques exploiting infinitesimals, recently analyzed in an article by
David Mumford, who wrote:
\begin{quote}
In my own education, I had assumed [that Enriques and the Italians]
were irrevocably stuck. \ldots{} As I see it now, Enriques must be
credited with a nearly complete geometric proof using, as did
Grothendieck, higher order infinitesimal deformations. \ldots{} Let's
be careful: he certainly had the correct ideas about infinitesimal
geometry, though he had no idea at all how to make precise definitions
(Mumford 2011 \cite{Mu11}).
\end{quote}

Another example is important work by Cauchy (see entry~\ref{ca} below)
on singular integrals and Fourier series using infinitesimals and
infinitesimally defined ``Dirac'' delta functions (these precede Dirac
by a century), which was forgotten for a number of decades because of
shifting foundational biases.  The presence of Dirac delta functions
in Cauchy's oeuvre was noted in (Freudenthal 1971 \cite{Fr71a}) and
analyzed by Laugwitz (1989 \cite{Lau89}), (1992a \cite{Lau92}); see
also (Katz \& Tall 2012 \cite{KT2}) and (Tall \& Katz 2013 \cite{TK}).

Recent papers on Leibniz (Katz \& Sherry \cite{KS1}, \cite{KS2};
Sherry \& Katz \cite{SK}) argue that, contrary to widespread
perceptions, Leibniz's system for infinitesimal calculus was
\emph{not} inconsistent (see entry~\ref{mr} on \emph{mathematical
rigor} for a discussion of the term).  The significance and coherence
of Berkeley's critique of infinitesimal calculus have been routinely
exaggerated.  Berkeley's sarcastic tirades against infinitesimals fit
well with the ontological limitations imposed by the A-approach
favored by many historians, even though Berkeley's opposition, on
empiricist grounds, to an infinitely divisible continuum is profoundly
at odds with the A-approach.

A recent study of Fermat (Katz, Schaps \& Shnider 2013 \cite{KSS13})
shows how the nature of his contribution to the calculus was distorted
in recent Fermat scholarship, similarly due to an ``evolutionary
dead-end" bias (see entry~\ref{35}).

The Marburg school of Hermann Cohen, Cassirer, Natorp, and others
explored the philosophical foundations of the infinitesimal method
underpinning the mathematized natural sciences.  Their versatile, and
insufficiently known, contribution is analyzed in (Mormann \& Katz
2013 \cite{MK}).

A number of recent articles have pioneered a re-evaluation of the
history and philosophy of mathematics, analyzing the shortcomings of
received views, and shedding new light on the deleterious effect of
the latter on the philosophy, the practice, and the applications of
mathematics.  Some of the conclusions of such a re-evaluation are
presented below.

\section
{Adequality to Chimeras}

Some topics from the history of infinitesimals illustrating our
approach appear below in alphabetical order.

\subsection{Adequality}
\label{ad}

Adequality is a technique used by Fermat to solve problems of
tangents, problems of maxima and minima, and other variational
problems.  The term \emph{adequality} derives from the \parisotes{} of
Diophantus (see entry~\ref{di}).  The technique involves an element of
approximation and ``smallness", represented by a small variation~$E$,
as in the familiar difference~$f(A+E)-f(A)$.  Fermat used adequality
in particular to find the tangents of transcendental curves like the
cycloid, that were considered to be ``mechanical" curves off-limits to
geometry, by Descartes.  Fermat also used it to solve the variational
problem of the refraction of light so as to obtain Snell's law.
Adequality incorporated a procedure of discarding higher-order terms
in~$E$ (\emph{without} setting them equal to zero).  Such a heuristic
procedure was ultimately formalized mathematically in terms of the
\emph{standard part principle} (entry~\ref{spf}) in A.~Robinson's
theory of infinitesimals starting with (Robinson 1961 \cite{Ro61}).
Fermat's adequality is comparable to Leibniz's \emph{transcendental
law of homogeneity}; see entry~\ref{tlh} and \cite{KSS13}.

\subsection{Archimedean axiom}
\label{aa}

What is known today as the Archimedean axiom first appears in Euclid's
\emph{Elements}, Book V, as Definition 4 (Euclid \cite{Euc},
definition~V.4).  It is exploited in (Euclid \cite{Euc},
Proposition~V.8).  We include bracketed symbolic notation so as to
clarify the definition:
\begin{quote}
Magnitudes [~$a,b$ ] are said to have a ratio with respect to one
another which, being multiplied [~$na$ ]
are capable of exceeding
one another [~$na>b$ ].
\end{quote}
It can be formalized as follows:%
\footnote{\label{f5}See e.g., the version of the Archimedean axiom in
(Hilbert 1899 \cite[p.~19]{Hi}).  Note that we have avoided using
``$0$'' in formula~\eqref{21c}, as in ``$\forall a>0$'', since~$0$ was
not part of the conceptual framework of the Greeks.  The term
``multiplied" in the English translation of Euclid's definition~V.4
corresponds to the Greek term
$\pi{}o\lambda\lambda\alpha\pi\lambda\alpha\sigma\iota
\alpha\zeta\acute{o}\mu\varepsilon\nu\alpha$.
%
%
A common formalisation of the noun ``multiple'',
$\pi{}o\lambda\lambda\alpha\pi\lambda\acute\alpha\sigma\iota{o}\nu$,
is~$na=a+...+a$.}
\begin{equation}
\label{21c}
(\forall a,b)(\exists n\in \N)\, [na>b],\ \ \mbox{where}\ \
na=\underbrace{a+...+a}_{n\;\text{times}}.
\end{equation}
Next, it appears in the papers of Archimedes as the following lemma
(see Archimedes \cite{Arc}, I, Lamb.~5):
\begin{quote}
Of unequal lines, unequal surfaces, and unequal solids [~$a,b,c$ ],
the greater exceeds
the lesser [~$a<b$ ] by such a magnitude [~$b-a$ ] as, when added to
itself [~$n(b-a)$~], can be made to exceed any assigned magnitude
[~$c$ ] among those which are comparable with one another (Heath 1897
\cite[p.~4]{He97}).
\end{quote}
This can be formalized as follows:
\begin{equation}
\label{22c}
(\forall{a,b,c})(\exists{n\in\N})\, [a<b\rightarrow n(b-a)>c].
\end{equation}
Note that Euclid's definition V.4 and the lemma of Archimedes are
\emph{not} logically equivalent (see entry~\ref{eu},
footnote~\ref{f5a}).

The Archimedean axiom plays no role in the plane geometry as developed
in Books I-IV of \emph{The Elements}.%
\footnote{With the exception of Proposition III.16 where so called
\emph{horn angles} appear.  These could be considered as
non-Archimedean magnitudes relative to rectilinear angles.}
Interpreting geometry in ordered fields, or in geometry over fields
for short, one knows that~$\mathbb{F}^2$ is a model of Euclid's plane,
where~$(\mathbb{F},+,\cdot,0,1,<)$ is a Euclidean field, i.e., an
ordered field closed under the square root operation.
Consequently,~$\R^*\times\R^*$ (where $\R^*$ is a hyperreal field) is
a model of Euclid's plane, as well (see entry~\ref{mo} on \emph{modern
implementations}).  Euclid's definition V.4 is discussed in more
detail in entry~\ref{eu}.

Otto Stolz rediscovered the Archimedean axiom for mathematicians,
making it one of his axioms for magnitudes and giving it the following
form: if~$a>b$, then there is a multiple of~$b$ such that~$nb>a$
(Stolz 1885 \cite[p.~69]{Sto}).%
\footnote{See Ehrlich \cite{Eh06} for additional historical details
concerning Stolz's account of the Archimedean Axiom.}
At the same time, in his development of the integers Stolz implicitly
used the Archimedean axiom.  Stolz's visionary realisation of the
importance of the Archimedean axiom, and his work on non-Archimedean
systems, stand in sharp contrast with Cantor's remarks on
infinitesimals (see entry~\ref{mr} on \emph{mathematical rigor}).

In modern mathematics, the theory of ordered fields employs the
following form of the Archimedean axiom (see e.g., Hilbert 1899
\cite[p.~27]{Hi}):
\[
(\forall x>0)\; (\forall \epsilon>0)\; (\exists n\in\N) \;
[n\epsilon>x],
\]
or equivalently
\begin{equation}
\label{21b} 
(\forall\epsilon>0)\;(\exists n\in\N)\;\;[n\epsilon>1].
\end{equation}
A number system satisfying~\eqref{21b} will be referred to as an
\emph{Archimedean continuum}.  In the contrary case, there is an
element~$\epsilon>0$ called an infinitesimal such that no finite sum
$\epsilon+\epsilon+\ldots+\epsilon$ will ever reach~$1$; in other
words,
\begin{equation}
\label{22b}
(\exists\epsilon>0)\;(\forall
n\in\N)\;\left[\epsilon\leq\tfrac{1}{n}\right].
\end{equation}
A number system satisfying~\eqref{22b} is referred to as a
\emph{Bernoullian continuum} (i.e., a non-Archimedean continuum); see
entry~\ref{be}.

\subsection{Berkeley, George}
\label{gb}

George Berkeley (1685-1753) was a cleric whose \emph{empiricist}
(i.e., based on sensations, or \emph{sensationalist}) metaphysics
tolerated no conceptual innovations, like infinitesimals, without an
\emph{empirical} counterpart or referent.  Berkeley was similarly
opposed, on metaphysical grounds, to infinite divisibility of the
continuum (which he referred to as \emph{extension}), an idea widely
taken for granted today.  In addition to his outdated
\emph{metaphysical criticism} of the infinitesimal calculus of Newton
and Leibniz, Berkeley also formulated a \emph{logical criticism}.%
\footnote{Berkeley's criticism was dissected into its logical and
metaphysical components in (Sherry 1987 \cite{She87}).}
Berkeley claimed to have detected a logical fallacy at the basis of
the method.  In terms of Fermat's~$E$ occuring in his
\emph{adequality} (entry~\ref{ad}), Berkeley's objection can be
formulated as follows:
\begin{quote}
The increment~$E$ is assumed to be nonzero at the beginning of the
calculation, but zero at its conclusion, an apparent logical fallacy.
\end{quote}
However,~$E$ is not \emph{assumed to be zero} at the end of the
calculation, but rather is \emph{discarded} at the end of the
calculation (see entry~\ref{24b} for more details).  Such a technique
was the content of Fermat's \emph{adequality} (see entry~\ref{ad}) and
Leibniz's \emph{transcendental law of homogeneity} (see
entry~\ref{tlh}), where the relation of equality has to be suitably
interpreted (see entry~\ref{re} on relation~$\adequal$).  The
technique is closely related to taking the \emph{limit} (of a typical
expression such as~$\frac{f(A+E)-f(A)}{E}$ for example) in
Weierstrass's approach, and to taking the \emph{standard part} (see
entry~\ref{spf}) in Robinson's approach.

Meanwhile, Berkeley's own attempt to explain the calculation of the
derivative of~$x^2$ in \emph{The Analyst} contains a logical
circularity.  Namely, Berkeley's argument relies on the determination
of the tangents of a parabola by Apollonius (which is eqivalent to the
calculation of the derivative).  This circularity in Berkeley's
argument was analyzed in (Andersen 2011 \cite{An11}).  


\subsection{Berkeley's logical criticism}
\label{24b}

Berkeley's \emph{logical criticism} of the calculus amounts to the
contention that the evanescent increment is first assumed to be
non-zero to set up an algebraic expression, and then \emph{treated as
zero} in \emph{discarding} the terms that contained that increment
when the increment is said to vanish.  In modern terms, Berkeley was
claiming that the calculus was based on an inconsistency of type
\[
(dx\not=0)\wedge(dx=0).
\]
The criticism, however, involves a misunderstanding of Leibniz's
method.  The rebuttal of Berkeley's criticism is that the evanescent
increment need \emph{not} be ``treated as zero'', but, rather, is
merely \emph{discarded} through an application of the
\emph{transcendental law of homogeneity} by Leibniz, as illustrated in
entry~\ref{pr} in the case of the product rule.

While consistent (in the sense of entry~\ref{mr}, level~\eqref{l2}),
Leibniz's system unquestionably relied on \emph{heuristic} principles
such as the laws of continuity and homogeneity, and thus fell short of
a standard of rigor \emph{if} measured by today's criteria (see
entry~\ref{mr} on \emph{mathematical rigor}).  On the other hand, the
consistency and resilience of Leibniz's system is confirmed through
the development of \emph{modern implementations} of Leibniz's
heuristic principles (see entry~\ref{mo}).

\subsection{Bernoulli, Johann}
\label{be}

Johann Bernoulli (1667-1748) was a disciple of Leibniz's who, having
learned an infinitesimal methodology for the calculus from the master,
never wavered from it.  This is in contrast to Leibniz himself, who,
throughout his career, used both
\begin{enumerate}
\item[(A)] an Archimedean methodology (proof by exhaustion), and
\item[(B)] an infinitesimal methodology,
\end{enumerate}
in a symbiotic fashion.  Thus, Leibniz relied on the A-methodology to
underwrite and justify the B-methodology, and he exploited the
B-methodology to shorten the path to discovery (\emph{Ars
Inveniendi}).  Historians often name Bernoulli as the first
mathematician to have adhered systematically to the infinitesimal
approach as the basis for the calculus.  We refer to an
infinitesimal-enriched number system as a \emph{B-continuum}, as
opposed to an Archimedean \emph{A-continuum}, i.e., a continuum
satisfying the \emph{Archimedean axiom} (see entry~\ref{aa}).

\subsection{Bishop, Errett}
\label{EB}

Errett Bishop (1928-1983) was a mathematical constructivist who,
unlike his fellow intuitionist%
\footnote{Bishop was not an Intuitionist in the narrow sense of the
term, in that he never worked with Brouwer's continuum or ``free
choice sequences".  We are using the term ``Intuitionism'' in a
broader sense (i.e., mathematics based on intuitionistic logic) that
incorporates constructivism, as used for example by Abraham Robinson
in the comment quoted at the end of this entry.}
\emph{Arend Heyting} (see entry~\ref{AH}), held a dim view of
classical mathematics in general and Robinson's infinitesimals in
particular.  Discouraged by the apparent non-constructivity of his
early work in functional analysis (notably \cite{Bi65A}), he believed
to have found the culprit in the law of excluded middle (LEM), the key
logical ingredient in every proof by contradiction.  He spent the
remaining 18 years of his life in an effort to expunge the reliance on
LEM (which he dubbed ``the principle of omniscience'' in \cite{Bi75})
from analysis, and sought to define \emph{meaning} itself in
mathematics in terms of such LEM-extirpation.

Accordingly, he described classical mathematics as both a
\emph{debasement of meaning} (Bishop 1973 \cite[p.~1]{Bi85}) and
\emph{sawdust} (Bishop 1973 \cite[p.~14]{Bi85}), and did not hesitate
to speak of both \emph{crisis} (Bishop 1975 \cite{Bi75}) and
\emph{schizophrenia} (Bishop 1973 \cite{Bi85}) in contemporary
mathematics, predicting an imminent demise of classical mathematics in
the following terms:
\begin{quote}
Very possibly classical mathematics will cease to exist as an
independent discipline (Bishop 1968 \cite[p.~54]{Bi68}).
\end{quote}
His attack in (Bishop 1977 \cite{Bi77}) on calculus pedagogy based on
Robinson's infinitesimals was a natural outgrowth of his general
opposition to the logical underpinnings of classical mathematics, as
analyzed in (Katz \& Katz 2011 \cite{KK11d}).  Robinson formulated a
brief but penetrating appraisal of Bishop's ventures into the history
and philosophy of mathematics as follows:
\begin{quote}
The sections of [Bishop's] book that attempt to describe the
philosophical and historical background of [the] remarkable endeavor
[of Intuitionism] are more vigorous than accurate and tend to belittle
or ignore the efforts of others who have worked in the same general
direction (Robinson 1968 \cite[p.~921]{Ro68}).
\end{quote}

See entry~\ref{AC} for a related criticism by Alain Connes.

\subsection{Cantor, Georg}

Georg Cantor (1845-1918) is familiar to the modern reader as the
\emph{underappreciated} creator of the ``Cantorian paradise" which
David Hilbert would not be expelled out of, as well as the tragic
hero, allegedly \emph{persecuted} by Kronecker, who ended his days in
a lunatic asylum.  Cantor historian J. Dauben notes, however, an
\emph{underappreciated} aspect of Cantor's scientific activity, namely
his principled \emph{persecution} of infinitesimalists:
\begin{quote}
Cantor devoted some of his most vituperative correspondence, as well
as a portion of the \emph{Beitr\"age}, to attacking what he described
at one point as the `infinitesimal Cholera bacillus of mathematics',
which had spread from Germany through the work of Thomae, du Bois
Reymond and Stolz, to infect Italian mathematics \ldots{} Any
acceptance of infinitesimals necessarily meant that his own theory of
number was incomplete.  Thus to accept the work of Thomae, du
Bois-Reymond, Stolz and Veronese was to deny the perfection of
Cantor's own creation.  Understandably, Cantor launched a thorough
campaign to discredit Veronese's work in every way possible (Dauben
1980 \cite[pp.~216-217]{Da80}).
\end{quote}
A discussion of Cantor's flawed investigation of the \emph{Archimedean
axiom} (see entry~\ref{aa}) may be found in entry~\ref{mr} on
\emph{mathematical rigor}.%
\footnote{\label{melancholy1}Cantor's dubious claim that the
infinitesimal leads to contradictions was endorsed by no less an
authority than B.~Russell; see footnote~\ref{f8} in entry~\ref{mr}.}

\subsection{Cauchy, Augustin-Louis}
\label{ca}

Augustin-Louis Cauchy (1789-1857) is often viewed in the history of
mathematics literature as a precursor of Weierstrass.  Note, however,
that contrary to a common misconception, Cauchy never gave an
$\epsilon, \delta$ definition of either limit or continuity (see
entry~\ref{vc} on \emph{variable quantity} for Cauchy's definition of
limit).  Rather, his approach to continuity was via what is known
today as \emph{microcontinuity} (see entry~\ref{co}).  Several recent
articles (B\l aszczyk et al.~\cite{BKS}; Borovik \& Katz \cite{BK};
Br\aa{}ting \cite{Br}; Katz \& Katz \cite{KK11b}, \cite{KK11a}; Katz
\& Tall \cite{KT2}; Tall \& Katz \cite{TK}), have argued that a
proto-Weierstrassian reading of Cauchy is one-sided and obscures
Cauchy's important contributions, including not only his infinitesimal
definition of continuity but also such innovations as his
infinitesimally defined (``Dirac'') delta function, with applications
in Fourier analysis and evaluation of singular integrals, and his
study of orders of growth of infinitesimals that anticipated the work
of Paul du Bois-Reymond, Borel, Hardy, and ultimately Skolem
(\cite{Sk33}, \cite{Sk34}, \cite{Sk55}) and Robinson.

To elaborate on Cauchy's ``Dirac'' delta function, note the following
formula from (Cauchy 1827 \cite[p.~188]{Ca27}) in terms of an
infinitesimal~$\alpha$:
\begin{equation}
\label{151}
\frac{1}{2} \int_{a-\epsilon}^{a+\epsilon} F(\mu) \frac{\alpha \;
d\mu}{\alpha^2 + (\mu-a)^2} = \frac{\pi}{2} F(a).
\end{equation}
Replacing Cauchy's expression~$\frac{\alpha}{\alpha^2 + (\mu-a)^2}$ by
$\delta_a(\mu)$, one obtains Dirac's formula up to trivial
modifications (see Dirac \cite[p.~59]{Dir}):
\[
\int_{-\infty}^{\infty} f(x)\delta(x)=f(0)
\]
Cauchy's 1853 paper on a notion closely related to uniform convergence
was recently examined in (Katz \& Katz 2011 \cite{KK11b}) and
(B\l{}aszczyk et al.~2012 \cite{BKS}).  Cauchy handles the said notion
using infinitesimals, including one generated by the null
sequence~$(\frac{1}{n})$.

Meanwhile, N\'u\~nez et al.~(1999 \cite[p.~54]{NEM}) coined the term
`Cauchy--Weierstrass definition of continuity'.  Since Cauchy gave an
infinitesimal definition and Weierstrass, an $\epsilon,\delta$ one,
such a coinage is an oxymoron.  J.~Gray (2008a \cite[p.~62]{Gr08a})
lists continuity among concepts Cauchy allegedly defined using
`limiting arguments', but Gray unfortunately confuses the term `limit'
as \emph{bound} with `limit' as in \emph{variable tending to a
quantity}, since the term `limits' appear in Cauchy's definition only
in the sense of endpoints (bounds) of an interval.  Not to be outdone,
Kline (1980 \cite[p. 273]{Kl80}) claims that ``Cauchy's work not only
banished [infinitesimals] but disposed of any need for them.''
Hawking (2007 \cite[p. 639]{Ha07}) does reproduce Cauchy's
infinitesimal definition, yet on the same page 639 claims that Cauchy
``was particularly concerned to banish infinitesimals," apparently
unaware of a comical non-sequitur he committed.

\subsection{Chimeras}
\label{AC}

Alain Connes (1947--) formulated criticisms of Robinson's
infinitesimals between the years 1995 and 2007, on at least seven
separate occasions (see Kanovei et al.~2012 \cite{KKM}, Section 3.1,
Table~1).  These range from pejorative epithets such as
``inadequate'', ``disappointing'', ``chimera'', and ``irremediable
defect'', to ``the end of the rope for being `explicit'\,''.

Connes sought to exploit the Solovay model~$\mathcal{S}$ (Solovay 1970
\cite{So}) as ammunition against non-standard analysis, but the model
tends to boomerang, undercutting Connes' own earlier work in
functional analysis.  Connes described the hyperreals as both a
``virtual theory" and a ``chimera'', yet acknowledged that his
argument relies on the transfer principle (see entry~\ref{mo}).
In~$\mathcal{S}$, all definable sets of reals are Lebesgue measurable,
suggesting that Connes views a theory as being ``virtual'' if it is
not \emph{definable} in a suitable model of ZFC.  If so, Connes' claim
that a theory of the hyperreals is ``virtual" is refuted by the
existence of a definable model of the hyperreal field (Kanovei \&
Shelah \cite{KS}).  Free ultrafilters aren't definable, yet Connes
exploited such ultrafilters both in his own earlier work on the
classification of factors in the 1970s and 80s, and in his magnum opus
\emph{Noncommutative Geometry} (Connes 1994 \cite[ch.~V,
sect.~6.$\delta$, Def.~11]{Co94}), raising the question whether the
latter may not be vulnerable to Connes' criticism of virtuality.  The
article \cite{KKM} analyzed the philosophical underpinnings of Connes'
argument based on G\"odel's incompleteness theorem, and detected an
apparent circularity in Connes' logic.  The article~\cite{KKM} also
documented the reliance on non-constructive foundational material, and
specifically on the Dixmier trace~$-\hskip-10pt\int$ (featured on the
front cover of Connes' \emph{magnum opus}) and the Hahn--Banach
theorem, in Connes' own framework; see also \cite{KL}.

See entry~\ref{EB} for a related criticism by Errett Bishop.


\section{Continuity to indivisibles}

\subsection{Continuity}
\label{co}

Of the two main definitions of continuity of a function, Definition A
(see below) is operative in either a B-continuum or an A-continuum
(satisfying the Archimedean axiom; see entry~\ref{aa}), while
Definition B only works in a B-continuum (i.e., an
infinitesimal-enriched or Bernoullian continuum; see entry~\ref{be}).
\begin{itemize}
\item
\emph{Definition A} ($\epsilon, \delta$ approach): A real function~$f$
is continuous at a real point~$x$ if and only if
\[
\quad\quad\quad(\forall\epsilon>0)\; (\exists\delta>0)\; (\forall x')
\;\big[\,|x-x'|<\delta \rightarrow |f(x)-f(x')|<\epsilon\,\big].
\]
\item
\emph{Definition B} (microcontinuity): A real function~$f$ is
continuous at a real point~$x$ if and only if
\begin{equation}
\label{26}
(\forall x')\; \big[\,x'\adequal x\,\rightarrow\,f(x')\adequal
f(x)\,\big].
\end{equation}
\end{itemize}
In formula~\eqref{26}, the natural extension of~$f$ is still
denoted~$f$, and the symbol~``$\adequal$'' stands for the relation of
being infinitely close; thus,~$x'\adequal x$ if and only if~$x'-x$ is
infinitesimal (see entry~\ref{re} on relation~``$\adequal$'').



\subsection{Diophantus}
\label{di}

Diophantus of Alexandria (who lived about 1800 years ago) contributed
indirectly to the development of infinitesimal calculus through the
technique called \parisotes{}, developed in his work
\emph{Arithmetica}, Book Five, problems 12, 14, and 17.  The term
\parisotes{} can be literally translated as ``approximate equality''.
This was rendered as \emph{adaequalitas} in Bachet's Latin translation
\cite{Bachet}, and \emph{ad\'egalit\'e} in French (see entry~\ref{ad}
on \emph{adequality}).  The term was used by Fermat to describe the
comparison of values of an algebraic expression, or what would today
be called a function~$f$, at nearby points~$A$ and~$A+E$, and to seek
extrema by a technique that would be re-formulated today in terms of
the vanishing of $\frac{f(A+E)-f(A)}{E}$ after discarding the
remaining~$E$-terms; see (Katz, Schaps \& Shnider 2013 \cite{KSS13}).

\subsection{Euclid's definition V.4}
\label{eu}

Euclid's definition V.4 was already mentioned in entry~\ref{aa}.  In
addition to Book V, it appears in Books X and XII and is used in the
method of exhaustion (see Euclid \cite{Euc}, Propositions X.1, XII.2).
The method of exhaustion was exploited intensively by both Archimedes
and Leibniz (see entry~\ref{38} on Leibniz's work \emph{De
Quadra\-tura}).  It was revived in the 19th century in the theory of
the Riemann integral.

Euclid's Book V sets the basis for the theory of similar figures
developed in Book VI.  Great mathematicians of the 17th century like
Descartes, Leibniz, and Newton exploited Euclid's theory of similar
figures of Book VI while paying no attention to its axiomatic
background.%
\footnote{\label{f3}Leibniz and Newton apparently applied Euclid's
conclusions in a context where the said conclusions did not
technically speaking apply, namely to infinitesimal figures such as
the characteristic triangle, i.e., triangle with sides~$dx$,~$dy$,
and~$ds$.}
Over time Euclid's Book~V became a subject of interest for historians
and editors
%
%
alone.

To formalize Definition V.4, one needed a formula for Euclid's notion
of ``multiple" and an idea of total order. Some progress in this
direction was made by Robert Simson in 1762.%
\footnote{See Simson's axioms that supplement the definitions of
Book~V as elaborated in (Simson 1762 \cite[p.~114-115]{Si}).}
In 1876, Hermann Hankel provided a modern reconstruction of Book~V.
Combining his own historical studies with an idea of order compatible
with addition developed by Hermann Grassmann (1861 \cite{Gras}), he
gave a formula that to this day is accepted as a formalisation of
Euclid's definition of proportion in V.5 (Hankel 1876
\cite[pp.~389-398]{Han}).  Euclid's proportion is a relation among
four ``magnitudes", such as
\[
A:B::C:D.
\]
It was interpreted by Hankel as the relation
\[
(\forall{m,n})\,\big[(nA\,{>}_1^{\phantom{I}}\,mB\rightarrow
nC>_2^{\phantom{I}}mD)\wedge
\]
\[
\wedge(nA=mB\rightarrow nC=mD) \wedge
(nA<_1^{\phantom{I}}mB\rightarrow nC<_2^{\phantom{I}}mD)\big],
\] 
where~$n, m$ are natural numbers.  The indices on the inequalities
emphasize the fact that the ``magnitudes"~$A,B$ have to be of ``the
same kind", e.g., line segments, whereas~$C,D$ could be of another
kind, e.g., triangles.

In 1880, J. L.~Heiberg in his edition of Archimedes' \emph{Opera
omnia}, in a comment on a lemma of Archimedes, cites Euclid's
definition~V.4, noting that these two are the same axioms (Heiberg
1880 \cite[p.~11]{Hei}).%
\footnote{\label{f5a}In point of fact, Euclid's axiom V.4 and
Archimedes' lemma are not equivalent from the logical viewpoint.
Thus, the additive semigroup of positive appreciable limited
hyperreals satisfies V.4 but not Archimedes' lemma.}
This is the reason why Euclid's definition~V.4 is commonly known as
the Archimedean axiom.  Today we formalize Euclid's definition~V.4 as
in \eqref{21c}, while the Archimedean lemma is rendered by
formula~\eqref{22c}.


\subsection{Euler, Leonhard}
\label{32}

Euler's \emph{Introductio in Analysin Infinitorum} (1748 \cite{ai})
contains remarkable calculations carried out in an extended number
system in which the basic algebraic operations are applied to
infinitely small and infinitely large quantities.  Thus, in Chapter 7,
``Exponentials and Logarithms Expressed through Series", we find a
derivation of the power series for~$a^z$ starting from the
formula~$a^\omega=1+k\omega$, for~$\omega$ infinitely small and then
raising the equation to the infinitely great power%
\footnote{Euler used the symbol~$i$ for the infinite power.  Blanton
replaced this by~$j$ in the English translation so as to avoid a
notational clash with the standard symbol for~$\sqrt{-1}$.}
$j=\frac{z}{\omega}$ for a finite (appreciable)~$z$ to give
\[
a^z= a^{j\omega}=(1+k\omega)^j
\]
and finally expanding the right hand side as a power series by means
of the binomial formula.  In the chapters following Euler finds
infinite product expansions factoring the power series expansion for
transcendental functions (see entry~\ref{34} for his infinite product
formula for sine).  By Chapter 10, he has the tools to sum the series
for~$\zeta(2)$ where~$\zeta(s)=\sum_n n^{-s}$.  He explicitly
calculates~$\zeta(2k)$ for~$k=1,\ldots,13$ as well as many other
related infinite series.

In Chapter 3 of his \emph{Institutiones Calculi Differentialis} (1755
\cite{ed}), Euler deals with the methodology of the calculus, such as
the nature of infinitesimal and infinitely large quantities.  We will
cite the English translation \cite{Eu55} of the Latin original
\cite{ed}.  Here Euler writes that
\begin{quote}
[e]ven if someone \emph{denies} that infinite numbers really exist in
this world, still in mathematical speculations there arise questions
to which answers cannot be given unless we admit an infinite number
(ibid., \S\,82) [emphasis added--the authors].
\end{quote}
Euler's approach, countenancing the possibility of \emph{denying} that
``infinite numbers really exist'', is consonant with a Leibnizian view
of infinitesimal and infinite quantities as ``useful fictions'' (see
Katz \& Sherry \cite{KS1}; Sherry \& Katz \cite{SK}).  Euler then
notes that ``an infinitely small quantity is nothing but a vanishing
quantity, and so it is really equal to 0'' (ibid., \S\,83).  The
``equality'' in question is an \emph{arithmetic} one (see below).

Similarly, Leibniz combined a view of infinitesimals as ``useful
fictions'' and \emph{inassignable quantities}, with a generalized
notion of ``equality'' which was an equality up to an incomparably
negligible term.  Leibniz sought to codify this idea in terms of his
\emph{transcendental law of homogeneity} (TLH); see entry~\ref{tlh}.
Thus, Euler's formulas like
\begin{equation}
\label{fo31}
a+dx=a
\end{equation}
(where~$a$ ``is any finite quantity''; ibid., \S\S\,86, 87) are
consonant with a Leibnizian tradition (cf.~formula~\eqref{adeq} in
entry~\ref{tlh}).  To explain formulas like~\eqref{fo31}, Euler
elaborated two distinct ways (arithmetic and geometric) of comparing
quantities in the following terms:

\begin{quote}
Since we are going to show that an infinitely small quantity is really
zero, we must meet the objection of why we do not always use the same
symbol 0 for infinitely small quantities, rather than some special
ones\ldots [S]ince we have two ways to compare them, either
\emph{arithmetic} or \emph{geometric}, let us look at the quotients of
quantities to be compared in order to see the difference.

If we accept the notation used in the analysis of the infinite,
then~$dx$ indicates a quantity that is infinitely small, so that both
$dx =0$ and~$a\,dx=0$, where~$a$ is any finite quantity.  Despite
this, the \emph{geometric} ratio~$a\,dx: dx$ is finite, namely~$a:1$.
For this reason, these two infinitely small quantities,~$dx$
and~$a\,dx$, both being equal to~$0$, cannot be confused when we
consider their ratio.  In a similar way, we will deal with infinitely
small quantities~$dx$ and~$dy$ (ibid., \S\,86, p. 51-52) [emphasis
added--the authors].
\end{quote}
Euler proceeds to clarify the difference between the arithmetic and
geometric comparisons as follows:
\begin{quote}
Let $a$ be a finite quantity and let $dx$ be infinitely small.
%
%
The arithmetic ratio of equals is clear: Since $ndx =0$, we have
\[
a \pm ndx - a = 0.
\]
On the other hand, the geometric ratio is clearly of equals, since
\begin{equation}
\label{32b}
\frac{a \pm ndx}{a} = 1.
\end{equation}
From this we obtain the well-known rule that \emph{the infinitely
small vanishes in comparison with the finite and hence can be
neglected} (Euler 1755 \cite[\S 87]{Eu55}) [emphasis in the
original--the authors].
\end{quote}
Like Leibniz, Euler considers more than one way of comparing
quantities.  Euler's formula~\eqref{32b} indicates that his geometric
comparison is procedurally identical with the Leibnizian TLH.  Namely,
Euler's geometric comparision of a pair of quantities amounts to their
ratio being infinitely close to $1$; the same is true for TLH.  Thus,
one has $a+dx=a$ in this sense for an appreciable~$a\not=0$, but
\emph{not} $dx=0$ (which is true only \emph{arithmetically} in Euler's
sense).  Euler's ``geometric'' comparison was dubbed ``the principle
of cancellation'' in (Ferraro 2004 \cite[p.~47]{Fe04}).

Euler proceeds to present the usual rules of infinitesimal calculus,
which go back to Leibniz, L'H\^opital, and the Bernoullis, such as
\begin{equation}
\label{fo32}
a\,dx^m + b\,dx^n=a\,dx^m
\end{equation}
provided~$m<n$ ``since~$dx^n$ vanishes compared with~$dx^m$'' (ibid.,
\S\,89), relying on his ``geometric'' equality.  Euler introduces a
distinction between infinitesimals of different order, and directly
\emph{computes}%
\footnote{\label{f7}Note that Euler does not ``prove that the
expression is equal to 1''; such \emph{indirect} proofs are a
trademark of the~$\epsilon, \delta$ approach.  Rather, Euler
\emph{directly} computes (what would today be formalized as the
standard part of) the expression, illustrating one of the advantages
of the B-methodology over the A-methodology.}
a ratio of the form
\[
\frac{dx\pm dx^2}{dx}=1\pm dx=1
\]
of two particular infinitesimals, assigning the value~$1$ to it
(ibid., \S\,88).  Euler concludes:
\begin{quote}
Although all of them [infinitely small quantities] are equal to 0,
still they must be carefully distinguished one from the other if we
are to pay attention to their mutual relationships, which has been
explained through a geometric ratio (ibid., \S\,89).
\end{quote}
The Eulerian hierarchy of orders of infinitesimals harks back to
Leibniz's work (see entry~\ref{nieu} on Nieuwentijt for a historical
dissenting view).  The remarkable lucidity of Euler's procedures for
dealing with infinitesimals has unfortunately not been appreciated by
all commentators.  Thus, J.~Gray interrupts his biography of Euler by
suddenly declaring: ``At some point it should be admitted that Euler's
attempts at explaining the foundations of calculus in terms of
differentials, which are and are not zero, are dreadfully weak'' (Gray
2008b \cite[p.~6]{Gr08b}) but provides no evidence whatsoever for his
dubious claim.

\subsection{Euler's infinite product formula for sine}
\label{34}

The fruitfulness of Euler's infinitesimal approach can be illustrated
by some of the remarkable applications he obtained.  Thus, Euler
derived an infinite product decomposition for the sine and sinh
functions of the following form:
\begin{eqnarray}
\label{eu1}
\sinh x&=&x\,
\left(1+\frac{x^2}{\pi^2}\right)
\,\left(1+\frac{x^2}{4\pi^2}\right)
\,\left(1+\frac{x^2}{9\pi^2}\right)
\,\left(1+\frac{x^2}{16\pi^2}\right)
\,\dots\,\quad\quad
\\[1ex]
\label{eu2}
\sin x&=&x\,
\left(1-\frac{x^2}{\pi^2}\right)
\,\left(1-\frac{x^2}{4\pi^2}\right)
\,\left(1-\frac{x^2}{9\pi^2}\right)
\,\left(1-\frac{x^2}{16\pi^2}\right)
\,\dots\,\quad\quad
\end{eqnarray}
Decomposition~\eqref{eu2} generalizes an infinite product formula for
$\frac{\pi}{2}$ due to Wallis \cite{Wa}.  Euler also summed the
inverse square
series:~$1+\frac14+\frac19+\frac1{16}+\ldots=\frac{\pi^2}6$ (see
\cite{mt}) and obtained additional identities.  A common feature of
these formulas is that Euler's computations involve not only
infinitesimals but also infinitely large natural numbers, which Euler
sometimes
treats as if they were ordinary natural numbers.%
\footnote{Euler's procedure is therefore consonant with the Leibnizian
\emph{law of continuity} (see entry~\ref{lc}) though apparently Euler
does not refer explicitly to the latter.}
Similarly, Euler treats infinite series as polynomials of a specific
infinite degree.

The derivation of \eqref{eu1} and \eqref{eu2} in (Euler 1748
\cite[\S\,156]{ai}) can be broken up into seven steps as
follows.
\vspace{1ex}

\emph{Step 1}.  Euler observes that
\begin{eqnarray}
\label{eu3}
2\sinh x &=& e^x-e^{-x} =
\left(1+\frac x{j}\right)^{j}-\left(1-\frac x{j}\right)^{j},
\end{eqnarray}
where~${j}$ (or ``$i$'' in Euler \cite{ai}) is an infinitely large
natural number.  To motivate the next step, note that the expression
$x^j-1=(x-1)(1+x+x^2+\ldots+x^{j-1})$ can be factored further as
$\prod_{k=0}^{j-1}(x-\zeta^k)$, where~$\zeta=e^{2\pi i/j}$; conjugate
factors can then be combined to yield a decomposition into real
quadratic terms as below.
\vspace{1ex} 

\emph{Step 2}.  Euler uses the fact that~$a^{j}-b^{j}$ is the product
of the factors
\begin{equation}
\label{eu4}
a^2+b^2-2ab\cos\frac{2k\pi}{j}\,,\quad\text{where}\quad k\ge1\,,
\end{equation}
together with the factor~$a-b$ and, if~${j}$ is an even number, the
factor~$a+b$, as well.
\vspace{1ex}

\emph{Step 3}.  Setting~$a=1+\frac x{j}$ and~$b=1-\frac x{j}$
in~\eqref{eu3}, Euler transforms expression~\eqref{eu4} into the form
\begin{equation}
\label{36b} 
2+2\frac{x^2}{{j}^2}-2\bigg(1-\frac{x^2}{{j}^2}\bigg)
\cos\frac{2k\pi}{j}\,.
\end{equation}

\emph{Step 4}.  Euler then replaces~\eqref{36b} by the expression
\begin{equation}
\label{37}
\frac{4k^2\pi^2}{{j}^2}
\bigg(1+\frac{x^2}{k^2\pi^2}-\frac{x^2}{{j}^2}\bigg)\,,
\end{equation}
justifying this step by means of the formula
\begin{equation}
\label{38b}
\cos\frac{2k\pi}{j} = 1-\frac{2k^2\pi^2}{{j}^2}.
\end{equation}

\emph{Step 5}.  Next, Euler argues that the difference~$e^x-e^{-x}$ is
divisible by the expression
\[
1+\frac{x^2}{k^2\pi^2}-\frac{x^2}{j^2}
\]
from \eqref{37}, where ``we omit the term~$\frac{x^2}{{j}^2}$ since
even when multiplied by~${j}$, it remains infinitely small'' (English
translation from \cite{aie}).
\vspace{1ex}

\emph{Step 6}.  As there is still a factor of~$a-b=2x/{j}$, Euler
obtains the final equality \eqref{eu1}, arguing that then ``the
resulting first term will be~$x$'' (in order to conform to the
Maclaurin series for~$\sinh x$).

\vspace{1ex}
\emph{Step 7}.  Finally, formula~\eqref{eu2} is obtained from
\eqref{eu1} by means of the substitution~$x\mapsto
ix$.\qed\vspace{1ex}

We will discuss modern formalisations of Euler's argument in
entry~\ref{e36}.

\subsection{Euler's sine factorisation formalized}
\label{e36}

Euler's argument in favor of \eqref{eu1} and \eqref{eu2} was
formalized in terms of a ``nonstandard'' proof in (Luxemburg 1973
\cite{lux}).  However, the formalisation in~\cite{lux} deviates from
Euler's argument beginning with steps 3 and~4, and thus circumvents
the most problematic steps 5 and~6.

A proof in the framework of modern nonstandard analysis, formalizing
Euler's argument step-by-step throughout, appeared in (Kanovei 1988
\cite{kan}); see also (McKinzie \& Tuckey 1997 \cite{mt}) and (Kanovei
\& Reeken 2004 \cite[Section 2.4a]{kr}).  This formalisation
interprets problematic details of Euler's argument on the basis of
general principles of modern nonstandard analysis, as well as general
analytic facts that were known in Euler's time.  Such principles and
facts behind some early proofs in infinitesimal calculus are sometimes
referred to as ``hidden lemmas'' in this context; see (Laugwitz
\cite{Lau87}, \cite{Lau89}), (McKinzie \& Tuckey 1997 \cite{mt}).
For instance, the ``hidden lemma'' behind Step 4 above is the fact
that for a fixed~$x$, the terms of the Maclaurin expansion of~$\cos x$
tend to~$0$ faster than a convergent geometric series, allowing one to
infer that the effect of the transformation of step 4 on the product
of the factors~\eqref{36b} is infinitesimal.  Some ``hidden lemmas" of
a different kind, related to basic principles of nonstandard analysis,
are discussed in \cite[pp~43ff.]{mt}.

What clearly stands out from Euler's argument is his explicit use of
infinitesimal expressions such as~\eqref{36b} and~\eqref{37}, as well
as the approximate formula~\eqref{38b} which holds ``up to'' an
infinitesimal of higher order.  Thus, Euler used infinitesimals
\emph{par excellence}, rather than merely ratios thereof, in a routine
fashion in some of his best work.

Euler's use of infinite integers and their associated infinite
products (such as the decomposition of the sine function) were
interpreted in Robinson's framework in terms of hyperfinite sets.
Thus, Euler's product of~$j$-infinitely many factors in~\eqref{eu2} is
interpreted as a hyperfinite product in \cite[formula~(9), p.~74]{kr}.
A hyperfinite formalisation of Euler's argument involving infinite
integers and their associated products illustrates the successful
remodeling of the \emph{arguments} (and not merely the \emph{results})
of classical infinitesimal mathematics, as discussed in
entry~\ref{mr}.


\subsection{Fermat, Pierre}
\label{35}

Pierre de Fermat (1601-1665) developed a pioneering technique known as
\emph{adequality} (see entry~\ref{ad}) for finding tangents to curves
and for solving problems of maxima and minima.  (Katz, Schaps \&
Shnider 2013 \cite{KSS13}) analyze some of the main approaches in the
literature to the method of adequality, as well as its source in the
\parisotes{} of Diophantus (see entry~\ref{di}).  At least some of the
manifestations of adequality, such as Fermat's treatment of
transcendental curves and Snell's law, amount to variational
techniques exploiting a small (alternatively, infinitesimal)
variation~$E$.  Fermat's treatment of geometric and physical
applications suggests that an aspect of approximation is inherent in
adequality, as well as an aspect of smallness on the part of~$E$.

Fermat's use of the term \emph{adequality} relied on Bachet's rendering
of Diophantus.  Diophantus coined the term \emph{parisotes} for
mathematical purposes.  Bachet performed a semantic calque in passing
from \emph{par-iso\=o} to \emph{ad-aequo}.  A historically significant
parallel is found in the similar role of, respectively, adequality and
the \emph{transcendental law of homogeneity} (see entry~\ref{tlh}) in
the work of, respectively, Fermat and Leibniz on the problems of
maxima and minima.

Breger (1994 \cite{Bre94}) denies that the idea of ``smallness" was
relied upon by Fermat.  However, a detailed analysis (see
\cite{KSS13}) of Fermat's treatment of the cycloid reveals that Fermat
did rely on issues of ``smallness" in his treatment of the cycloid,
and reveals that Breger's interpretation thereof contains both
mathematical errors and errors of textual analysis.  Similarly,
Fermat's proof of Snell's law, a variational principle, unmistakably
relies on ideas of ``smallness".

Cifoletti (1990 \cite{Ci}) finds similarities between Fermat's
adequality and some procedures used in smooth infinitesimal analysis
of Lawvere and others.  Meanwhile, (J. Bell 2009 \cite{Bel09}) seeks
the historical sources of Lawvere's infinitesimals mainly in
Nieuwentijt (see entry~\ref{nieu}).





\subsection{Heyting, Arend}
\label{AH}
Arend Heyting (1898-1980) was a mathematical Intuitionist whose
lasting contribution was the formalisation of the Intuitionistic logic
underpinning the Intuitionism of his teacher Brouwer.  While Heyting
never worked on any theory of infinitesimals, he had several
opportunities to present an expert opinion on Robinson's theory.
Thus, in 1961, Robinson made public his new idea of non-standard
models for analysis, and ``communicated this almost immediately to
\ldots Heyting'' (see Dauben \cite[p.~259]{Da03}).  Robinson's first
paper on the subject was subsequently published in {\em Proceedings of
the Netherlands Royal Academy of Sciences\/}~\cite{Ro61}.  Heyting
praised non-standard analysis as ``a standard model of important
mathematical research'' (Heyting 1973 \cite[p.~136]{He73}).
Addressing Robinson, he declared:
\begin{quote}
you connected this extremely abstract part of model theory with a
theory apparently so far apart as the elementary calculus.  In doing
so you threw new light on the history of the calculus by giving a
clear sense to Leibniz's notion of infinitesimals (ibid).
\end{quote}
Intuitionist Heyting's admiration for the application of Robinson's
infinitesimals to calculus pedagogy is in stark contrast with the
views of his fellow constructivist E.~Bishop (entry~\ref{EB}).



\subsection{Indivisibles versus Infinitesimals}
\label{36}

Commentators use the term \emph{infinitesimal} to refer to a variety
of conceptions of the infinitely small, but the variety is not always
acknowledged.  It is important to distinguish the infinitesimal
methods of Archimedes and Cavalieri from those employed by Leibniz and
his successors.  To emphasize this distinction, we will say that
tradition prior to Leibniz employed \emph{indivisibles}.  For example,
in his heuristic proof that the area of a parabolic segment is 4/3 the
area of the inscribed triangle with the same base and vertex,
Archimedes imagines both figures to consist of perpendiculars of
various heights erected on the base.  The perpendiculars are
indivisibles in the sense that they are limits of division and so one
dimension less than the area.  In the same sense, the indivisibles of
which a line consists are points, and the indivisibles of which a
solid consists are planes.

Leibniz's infinitesimals are not indivisibles, for they have the same
dimension as the figures they form.  Thus, he treats curves as
composed of infinitesimal line intervals rather than indivisible
points.  The strategy of treating infinitesimals as dimensionally
homogeneous with the objects they compose seems to have originated
with Roberval or Torricelli, Cavalieri's student, and to have been
explicitly arithmetized in (Wallis 1656 \cite{Wa}).

Zeno's \emph{paradox of extension} admits of resolution in the
framework of Leibnizian infinitesimals (see entry~\ref{ze}).
Furthermore, only with the dimensionality retained is it possible to
make sense of the fundamental theorem of calculus, where one must
think about the rate of change of the \emph{area} under a curve,
another reason why indivisibles had to be abandoned in favor of
infinitesimals so as to enable the development of the calculus (see
Ely 2012 \cite{El12}).



\section{Leibniz to Nieuwentijt}

\subsection{Leibniz, Gottfried}

Gottfried Wilhelm Leibniz (1646-1716), the co-inventor of
infinitesimal calculus, is a key player in the parallel infinitesimal
track referred to by Felix Klein \cite[p.~214]{Kl08} (see
Section~\ref{one}).

Leibniz's \emph{law of continuity} (see entry~\ref{lc}) together with
his \emph{transcendental law of homogeneity} (which he already
discussed in his response to Nieuwentijt in 1695 as noted by
M.~Parmentier~\cite[p.~38]{Le89}, and later in greater detail in a
1710 article \cite{Le10b} cited in the seminal study of Leibnizian
methodology by H. Bos \cite{Bos}) form a basis for implementing the
calculus in the context of a B-continuum.

Many historians of the calculus deny significant continuity between
infinitesimal calculus of the 17th century and 20th century
developments such as Robinson's theory (see further discussion in Katz
\& Sherry \cite{KS1}).  Robinson's hyperreals require the resources of
modern logic; thus many commentators are comfortable denying a
historical continuity.  A notable exception is Robinson himself, whose
identification with the Leibnizian tradition inspired Lakatos,
Laugwitz, and others to consider the history of the infinitesimal in a
more favorable light.  Many historians have overestimated the force of
Berkeley's criticisms (see entry~\ref{gb}), by underestimating the
mathematical and philosophical resources available to Leibniz.

Leibniz's infinitesimals are fictions, not logical fictions, as
(Ishiguro 1990 \cite{Is}) proposed, but rather pure fictions, like
imaginaries, which are not eliminable by some syncategorematic
paraphrase; see (Sherry \& Katz 2013 \cite{SK}) and entry~\ref{38}
below.

In fact, Leibniz's defense of infinitesimals is more firmly grounded
than Berkeley's criticism thereof.  Moreover, Leibniz's system for
differential calculus was free of logical fallacies (see
entry~\ref{24b}).  This strengthens the conception of modern
infinitesimals as a formalisation of Leibniz's strategy of relating
inassignable to assignable quantities by means of his
\emph{transcendental law of homogeneity} (see entry~\ref{tlh}).

\subsection{Leibniz's \emph{De Quadratura}}
\label{38}

In 1675 Leibniz wrote a treatise on his infinitesimal methods,
\emph{On the Arithmetical Quadrature of the Circle, the Ellipse, and
the Hyperbola}, or \emph{De Quadratura}, as it is widely known.
However, the treatise appeared in print only in 1993 in a text edited
by Knobloch (Leibniz \cite{Le1993}).

\emph{De Quadratura} was interpreted by R. Arthur \cite{Art} and
others as supporting the thesis that Leibniz's infinitesimals are mere
shortcuts, eliminable by long-winded paraphrase.  This so-called
syncategorematic interpretation of Leibniz's calculus has gained a
number of adherents.  We believe this interpretation to be in error.
In the first place, Leibniz wrote the treatise at a time when
infinitesimals were despised by the French Academy, a society whose
approval and acceptance he eagerly sought.  More importantly, as
(Jesseph 2013 \cite{Je11}) has shown, \emph{De Quadratura} depends on
infinitesimal resources in order to construct an approximation to a
given curvilinear area less than any previously specified error.  This
problem is reminiscent of the difficulty that led to infinitesimal
methods in the first place.  Archimedes' method of exhaustion required
one to determine a value for the quadrature \emph{in advance} of
showing, by \emph{reductio} argument, that any departure from that
value entails a contradiction.  Archimedes possessed a heuristic,
indivisible method for finding such values, and the results were
justified by exhaustion, but only after the fact.  By the same token,
the use of infinitesimals is `just' a shortcut only if it is entirely
eliminable from quadratures, tangent constructions, etc.  Jesseph's
insight is that this is not the case.

Finally, the syncategorematic interpretation misrepresents a crucial
aspect of Leibniz's mathematical philosophy.  His conception of
mathematical fiction includes imaginary numbers, and he often sought
approbation for his infinitesimals by comparing them to imaginaries,
which were largely uncontroversial.  There is no suggestion by Leibniz
that imaginaries are eliminable by long-winded paraphrase.  Rather, he
praises imaginaries for their capacity to achieve universal harmony by
the greatest possible systematisation, and this characteristic is more
central to Leibniz's conception of infinitesimals than the idea that
they are mere shorthand. Just as imaginary roots both unified and
extended the method for solving cubics, likewise infinitesimals
unified and extended the method for quadrature so that, e.g.,
quadratures of general parabolas and hyperbolas could be found by the
same method used for quadratures of less difficult curves.  See also
(Tho 2012 \cite{Th}).


\subsection
{\emph{Lex continuitatis}}
\label{lc}

A heuristic principle called \emph{The law of continuity} (LC) was
formulated by Leibniz and is a key to appreciating Leibniz's vision of
infinitesimal calculus.  The LC asserts that whatever succeeds in the
finite, succeeds also in the infinite.  This form of the principle
appeared in a letter to Varignon (Leibniz 1702 \cite{Le02}).  A more
detailed form of LC in terms of the concept of \emph{terminus}
appeared in his text \emph{Cum Prodiisset}:
\begin{quote}
In any supposed continuous transition, ending in any terminus, it is
permissible to institute a general reasoning, in which the final
terminus may also be included (Leibniz 1701 \cite[p.~40]{Le01c})
\end{quote}

To elaborate, the LC postulates that whatever properties are satisfied
by ordinary or assignable quantities, should also be satisfied by
inassignable quantities (see entry~\ref{vc}) such as infinitesimals
(see Figure~\ref{LCLH}).  Thus, the trigonometric
formula~$\sin^2x+\cos^2x=1$ should be satisfied for an inassignable
(e.g., infinitesimal) input~$x$, as well.  In the 20th century this
heuristic principle was formalized as the \emph{transfer principle}
(see entry~\ref{mo}) of \Los{}--Robinson.

The significance of LC can be illustated by the fact that a failure to
take note of the law of continuity often led scholars astray.  Thus,
Nieuwentijt (see entry~\ref{nieu}) was led into something of a
dead-end with his nilpotent infinitesimals (ruled out by LC) of the
form~$\frac{1}{\infty}$.  J.~Bell's view of Nieuwentijt's approach as
a precursor of nilsquare infinitesimals of Lawvere (see Bell 2009
\cite{Bel09}) is plausible, though it could be noted that Lawvere's
nilsquare infinitesimals cannot be of the form~$\frac{1}{\infty}$.

\begin{figure}
\begin{equation*}
\left\{\begin{matrix}assignable \cr quantities\end{matrix}\right\}
\buildrel{\rm LC}\over\leadsto
\left\{\begin{matrix}inassignable \cr quantities\end{matrix}\right\}
\buildrel{\rm TLH}\over\leadsto
\left\{\begin{matrix}assignable \cr quantities\end{matrix}\right\}
\end{equation*}
\caption{\textsf{Leibniz's law of continuity (LC) takes one from
assignable to inassignable quantities, while his transcendental law of
homogeneity (TLH; entry~\ref{tlh}) returns one to assignable
quantities.}}
\label{LCLH}
\end{figure}

\subsection
{\emph{Lex homogeneorum transcendentalis}}
\label{tlh}

Leibniz's \emph{transcendental law of homogeneity}, or \emph{lex
homogeneorum transcendentalis} in the original Latin (Leibniz 1710
\cite{Le10b}), governs equations involving differentials.  Leibniz
historian H.~Bos interprets it as follows:
\begin{quote}
A quantity which is infinitely small with respect to another quantity
can be neglected if compared with that quantity.  Thus all terms in an
equation except those of the highest order of infinity, or the lowest
order of infinite smallness, can be discarded.  For instance,
\begin{equation}
\label{adeq}
a+dx =a
\end{equation}
\[
dx+ddy=dx
\]
etc.  The resulting equations satisfy this \dots requirement of
homogeneity (Bos 1974 \cite[p.~33]{Bos}).
\end{quote}
For an interpretation of the equality sign in the formulas above, see
entry~\ref{re} on the relation~$\adequal$.

The TLH associates to an inassignable quantity (such as~$a+dx$), an
assignable one (such as~$a$); see Figure~\ref{LCLH} for a relation
between LC and TLH.



\subsection{Mathematical rigor}
\label{mr}


There is a certain lack of clarity in the historical literature with
regard to issues of fruitfulness, consistency, and rigorousness of
mathematical writing.  As a rough guide, and so as to be able to
formulate useful distinctions when it comes to evaluating mathematical
writing from centuries past, we would like to consider three levels of
judging mathematical writing:
\begin{enumerate}
\item potentially fruitful but (logically) inconsistent;
\item
\label{l2}
(potentially) consistent but informal;
\item
formally consistent and fully rigorous according to currently
prevailing standards.
\end{enumerate}
As an example of level~(1) we would cite the work of Nieuwentijt
(entry~\ref{nieu}; see there for a discussion of the inconsistency).
Our prime example of level~(2) is provided by the Leibnizian
\emph{laws of continuity and homogeneity} (entries~\ref{lc} and
\ref{tlh}), which found rigorous implementation at level~(3) only
centuries later (see entry~\ref{mo} on \emph{modern implementations}).

A foundation rock of the received history of mathematical analysis is
the belief that mathematical rigor emerged starting in the 1870s
through the efforts of Cantor, Dedekind, Weierstrass, and others,
thereby replacing formerly unrigorous work of infinitesimalists from
Leibniz onward.  The philosophical underpinnings of such a belief were
analyzed in (Katz \& Katz 2012a \cite{KK11a}) where it was pointed out
that in mathematics, as in other sciences, former errors are
eliminated through a process of improved conceptual understanding,
evolving over time, of the key issues involved in that science.  

Thus, no scientific development can be claimed to have attained
perfect clarity or rigor merely on the grounds of having eliminated
earlier errors.  Moreover, no such claim for a single scientific
development is made either by the practitioners or by the historians
of the natural sciences.  It was further pointed out in~\cite{KK11a}
that the term \emph{mathematical rigor} itself is ambiguous, its
meaning varying according to context.  Four possible meanings for the
term were proposed in \cite{KK11a}:

\begin{enumerate}
\item
it is a shibboleth that identifies the speaker as belonging to a
clan of professional mathematicians;
\item
it represents the idea that as a scientific field develops, its
practitioners attain greater and more conceptual understanding of key
issues, and are less prone to error;
\item
it represents the idea that a search for greater correctness in
analysis \emph{inevitably} led Weierstrass specifically to epsilontics
(i.e., the A-approach) in the 1870s;
\item
it refers to the establishment of what are perceived to be the
ultimate foundations for mathematics by Cantor, eventually explicitly
expressed in axiomatic form by Zermelo and Fraenkel.
\end{enumerate}

Item~(1) may be pursued by a fashionable academic in the social
sciences, but does not get to the bottom of the issue.  Meanwhile,
item~(2) would be agreed upon by historians of the other sciences.

In this context, it is interesting to compare the investigation of the
Archi\-me\-dean property as performed by the would-be rigorist Cantor,
on the one hand, and the infinitesimalist Stolz, on the other.  Cantor
sought to derive the Archimedean property as a consequence of those of
a linear continuum.  Cantor's work in this area was not only
unrigorous but actually erroneous, whereas Stolz's work was fully
rigorous and even visionary.  Namely, Cantor's arguments ``proving''
the inconsistency of infinitesimals were based on an implicit
assumption of what is known today as the Kerry-Cantor axiom (see
Proietti 2008 \cite{Pr}).  Meanwhile, Stolz was the first modern
mathematician to realize the importance of the \emph{Archimedean
axiom} (see entry~\ref{aa}) as a separate axiom in its own right (see
Ehrlich 2006~\cite{Eh06}), and moreover developed some non-Archimedean
systems (Stolz 1885 \cite{Sto}).

In his \emph{Grundlagen der Geometrie} (Hilbert 1899 \cite{Hi}),
Hilbert did not develop a new geometry, but rather remodeled Euclid's
geometry.  More specifically Hilbert brought rigor into Euclid's
geometry, in the sense of formalizing both Euclid's
\emph{propositions} and Euclid's style of \emph{procedures and style
of reasoning}.

Note that Hilbert's system works for geometries built over a
non-Archimedean field, as Hilbert was fully aware.  Hilbert (1900
\cite[p.~207]{Hi00}) cites Dehn's counterexamples to Legendre's
theorem in the absence of the Archimedean axiom.  Dehn planes built
over a non-Archimedean field were used to prove certain cases of the
independence of Hilbert's axioms%
\footnote{\label{f8}It is a melancholy comment to note that, fully
three years later, the philosopher-mathematician Bertrand Russell was
still claiming, on Cantor's authority, that the infinitesimal ``leads
to contradictions'' (Russell 2003 \cite[p.~345]{Ru03}).  This set the
stage for several decades of anti-infinitesimal vitriol, including the
\emph{saline solution} of Parkhurst and Kingsland (see
Section~\ref{one}.)}
(see Cerroni 2007 \cite{Ce07}).

Robinson's theory similarly formalized 17th and 18th century analysis
by remodeling both its \emph{propositions} and its \emph{procedures
and reasoning}.  Using Weierstrassian~$\epsilon, \delta$ techniques,
one can recover only the \emph{propositions} but not the proof
procedures.  Thus, Euler's \emph{result} giving an infinite product
formula for sine (entry~\ref{34}) admits of numerous proofs in a
Weierstrassian context, but Robinson's framework provides a suitable
context in which Euler's \emph{proof}, relying on infinite integers,
can also be recovered.  This is the crux of the historical debate
concerning~$\epsilon, \delta$ \emph{versus} infinitesimals.  In short,
Robinson did for Leibniz what Hilbert did for Euclid.  Meanwhile,
epsilontists failed to do for Leibniz what Robinson did for Leibniz,
namely formalizing the procedures and reasoning of the historical
infinitesimal calculus.  This theme is pursued further in terms of the
internal/external distinction in entry~\ref{vc} on \emph{variable
quantity}.

\subsection{Modern implementations}
\label{mo}

In the 1940s, Hewitt~\cite{Hew48} developed a modern implementation of
an infinitesimal-enriched continuum extending~$\R$, by means of a
technique referred to today as the ultrapower construction.  We will
denote such an infinitesimal-enriched continuum by the new
symbol~$\RRR$ (``thick-R'').%
\footnote{A more traditional symbol is~${}^*\R$ or~$\R^*$.}
In the next decade, (\Los{} 1955 \cite{Lo}) proved his celebrated
theorem on ultraproducts, implying in particular that elementary (more
generally, first-order) statements over~$\R$ are true if and only if
they are true over~$\RRR$, yielding a modern implementation of the
Leibnizian \emph{law of continuity} (entry~\ref{lc}).  \Los's theorem
is equivalent to what is known in the literature as the \emph{transfer
principle}; see Keisler~\cite{Ke08}.  Every finite element of~$\RRR$
is infinitely close to a unique real number; see entry~\ref{spf} on
the \emph{standard part principle}.  Such a principle is a
mathematical implementation of Fermat's \emph{adequality}
(entry~\ref{ad}); of Leibniz's \emph{transcendental law of
homogeneity} (see entry~\ref{tlh}); and of Euler's \emph{principle of
cancellation}
(see discussion between formulas~\eqref{fo31} and \eqref{fo32} in
entry~\ref{32}).

\subsection{Nieuwentijt, Bernard}
\label{nieu}

Nieuwentijt%
\footnote{Alternative spellings are Nieuwentijdt or Nieuwentyt.}
(1654-1718) wrote a book \emph{Analysis Infinitorum} (1695) proposing
a system containing an infinite number, as well as infinitesimal
quantities formed by dividing finite numbers by this infinite one.
Nieuwentijt postulated that the product of two infinitesimals should
be exactly equal to zero.  In particular, an infinitesimal quantity is
nilpotent.  In an exchange of publications with Nieuwentijt on
infinitesimals (see Mancosu 1996 \cite[p.~161]{Ma96}), Leibniz and
Hermann claimed that this system is consistent only if all
infinitesimals are equal, rendering differential calculus useless.
Leibniz instead advocated a system in which the product of two
infinitesimals is incomparably smaller than either infinitesimal.
Nieuwentijt's objections compelled Leibniz in 1696 to elaborate on the
hierarchy of infinite and infinitesimal numbers entailed in a robust
infinitesimal system.

Nieuwentijt's nilpotent infinitesimals of the form~$\frac{1}{\infty}$
are ruled out by Leibniz's \emph{law of continuity} (entry~\ref{lc}).
J.~Bell's view of Nieuwentijt's approach as a precursor of nilsquare
infinitesimals of Lawvere (see Bell 2009 \cite{Bel09}) is plausible,
though it could be noted that Lawvere's nilsquare infinitesimals
cannot be of the form~$\frac{1}{\infty}$.



\section{Product rule to Zeno}

\subsection{Product rule}
\label{pr}

In the area of Leibniz scholarship, the received view is that
Leibniz's infinitesimal system was logically faulty and contained
internal contradictions allegedly exposed by the cleric \emph{George
Berkeley} (entry~\ref{gb}).  Such a view is fully compatible with the
A-track-dominated outlook, bestowing supremacy upon the reconstruction
of analysis accomplished through the efforts of Cantor, Dedekind,
Weierstrass, and their rigorous followers (see entry~\ref{mr} on
\emph{mathematical rigor}).  Does such a view represent an accurate
appraisal of Leibniz's system?

The articles (Katz \& Sherry 2012 \cite{KS2}; 2013 \cite{KS1}; Sherry
\& Katz \cite{SK}) building on the earlier work (Sherry 1987
\cite{She87}), argued that Leibniz's system was in fact consistent (in
the sense of level~\eqref{l2} of entry~\ref{mr}),%
\footnote{Concerning the status of Leibniz's system for differential
calculus, it may be more accurate to assert that it was \emph{not
inconsistent}, in the sense that the contradictions alleged by
Berkeley and others turn out not to have been there in the first
place, once one takes into account Leibniz's generalized notion of
equality and his transcendental law of homogeneity.}
and featured resilient heuristic principles such as the \emph{law of
continuity} (entry~\ref{lc}) and the \emph{transcendental law of
homogeneity} (TLH) (entry~\ref{tlh}), which were implemented in the
fullness of time as precise mathematical principles guiding the
behavior of modern infinitesimals.

How did Leibniz exploit the TLH in developing the calculus?  We will
now illustrate an application of the TLH in the particular example of
the derivation of the product rule.  The issue is the justification of
the last step in the following calculation:
\begin{equation}
\label{41}
\begin{aligned}
d(uv) &= (u+du)(v+dv)-uv=udv+vdu+du\,dv \\ & =udv+vdu.
\end{aligned}
\end{equation}

The last step in the calculation~\eqref{41}, namely
\[
{udv+vdu} + {du\,dv} = {udv+vdu}
\]
is an application of the TLH.%
\footnote{Leibniz had two laws of homogeneity, one for dimension and
the other for the order of infinitesimalness.  Bos notes that they
`disappeared from later developments' \cite[p.~35]{Bos}, referring to
Euler and Lagrange.  Note, however, the similarity to Euler's
\emph{principle of cancellation} (see Bair et al.~\cite{B11}).}

In his 1701 text \emph{Cum Prodiisset} \cite[p.~46-47]{Le01c}, Leibniz
presents an alternative justification of the product rule (see Bos
\cite[p.~58]{Bos}).  Here he divides by~$dx$ and argues with
differential quotients rather than differentials.  Adjusting Leibniz's
notation to fit with~\eqref{41}, we obtain an equivalent calculation%
\footnote{The special case treated by Leibniz is~$u(x)=x$.  This
limitation does not affect the conceptual structure of the argument.}
\[
\begin{aligned}[l]
\frac{d(uv)}{dx} &= \frac{(u+du)(v+dv)-uv}{dx} =
\frac{udv+vdu+du\,dv}{dx} \\ &= \frac{udv+vdu}{dx} + \frac{du\,dv}{dx}
= \frac{udv+vdu}{dx}
\end{aligned}
\]
Under suitable conditions the term~$\left(\frac{du\,dv}{dx}\right)$ is
infinitesimal, and therefore the last step
\begin{equation}
\label{udv}
\frac{udv+vdu}{dx}+\frac{du\,dv}{dx}=u\,\frac{dv}{dx}+v\,\frac{du}{dx}
\end{equation}
is legitimized as a special case of the TLH.  The TLH interprets the
equality sign in~\eqref{udv} and~\eqref{adeq} as the relation of being
infinitely close, i.e., an equality up to infinitesimal error.

\subsection{Relation~$\adequal$}
\label{re}

\begin{figure}
\[
\xymatrix@C=95pt{{} \ar@{-}[rr] \ar@{-}@<-0.5pt>[rr]
\ar@{-}@<0.5pt>[rr] & {} \ar@{->>}[d]^{\hbox{st}} & \hbox{\quad
B-continuum} \\ {} \ar@{-}[rr] & {} & \hbox{\quad A-continuum} }
\]
\caption{\textsf{Thick-to-thin: applying the law of homogeneity, or
taking standard part (the thickness of the top line is merely
conventional)}}
\label{31}
\end{figure}

Leibniz did not use our equality symbol but rather the symbol
``$\adequal$'' (McClenon 1923 \cite[p.~371]{Mc23}).  Using such a
symbol to denote the relation of being infinitely close, one could
write the calculation of the derivative of~$y=f(x)$ where~$f(x)=x^2$
as follows:
\[
\begin{aligned}
f'(x)&\adequal\frac{dy}{dx}
\\&
=\frac{(x+dx)^2-x^2}{dx}
\\&
=\frac{(x+dx+x)(x+dx-x)}{dx}
\\&
=2x+dx
\\&
\adequal 2x.
\end{aligned}
\]
Such a relation is formalized by the \emph{standard part function};
see entry~\ref{spf} and Figure~\ref{31}.



\subsection{Standard part principle}
\label{spf}

In any totally ordered field extension~$E$ of~$\R$, every finite
element~$x\in E$ is infinitely close to a suitable unique element
$x_0\in\R$.  Indeed, via the total order, the element~$x$ defines a
Dedekind cut on~$\R$, and the cut specifies a real
number~$x_0\in\R\subset E$.  The number $x_0$ is infinitely close
to~$x\in E$.  The subring~$E_f\subset E$ consisting of the finite
elements of~$E$ therefore admits a map
\[
\text{st}:E_f\to\R,\;x\mapsto{}x_0,
\]
called the \emph{standard part function}.

\begin{figure}
\includegraphics[height=1.7in]{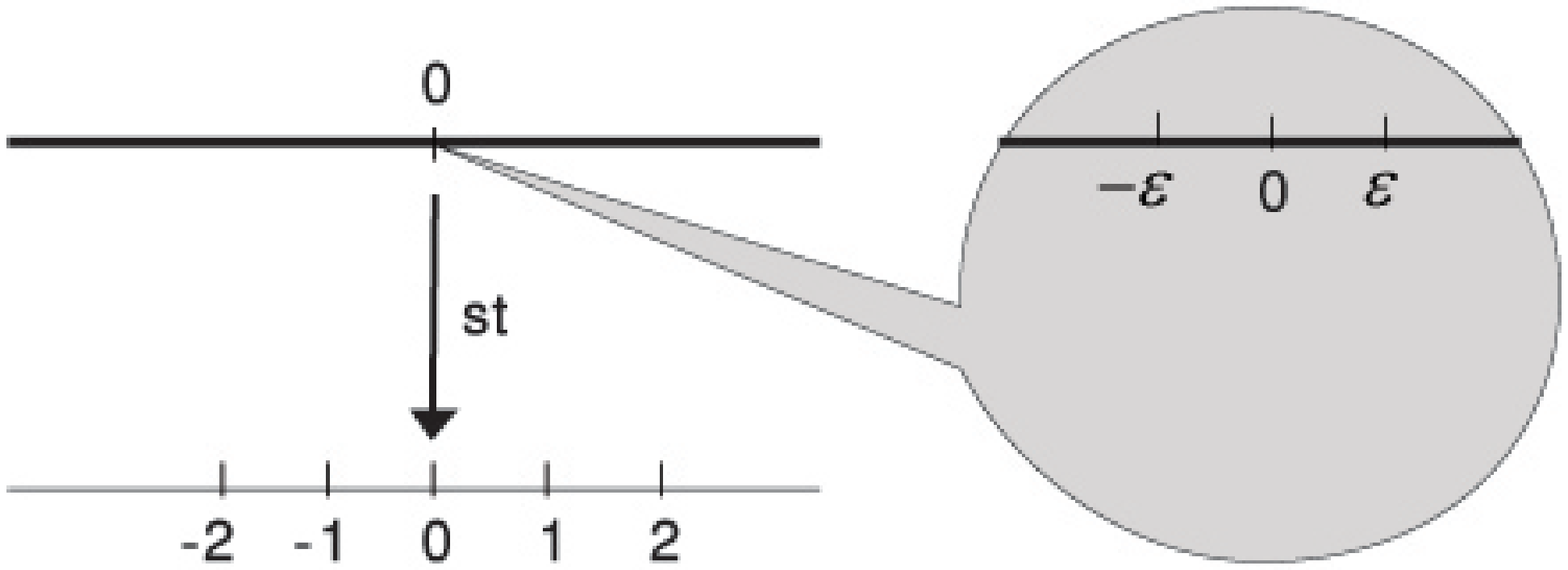}
\caption{\textsf{Zooming in on infinitesimal~$\varepsilon$ (here st$(\pm
\varepsilon)=0$).  The standard part function associates to every finite
hyperreal, the unique real number infinitely close to it. The bottom
line represents the ``thin" real continuum. The line at top represents
the ``thick" hyperreal continuum. The ``infinitesimal microscope" is
used to view an infinitesimal neighborhood of~$0$.  The derivative
$f'(x)$ of~$f(x)$ is then defined by the relation~$f'(x) \adequal
\frac{f(x+\varepsilon)-f(x)}{\varepsilon}$.}}
\label{tamar}
\end{figure}

The standard part function is illustrated in Figure~\ref{31}.  A more
detailed graphic representation may be found in Figure~\ref{tamar}.%
\footnote{For a recent study of optical diagrams in non-standard
analysis, see (Dossena \& Magnani \cite{DM}, \cite{MD}) and (Bair \&
Henry \cite{BH08}).}

The key remark, due to Robinson, is that the limit in the A-approach
and the standard part function in the B-approach are essentially
equivalent tools.  More specifically, the limit of a Cauchy
sequence~$(u_n)$ can be expressed, in the context of a hyperreal
enlargement of the number system, as the standard part of the
value~$u_H$ of the natural extension of the sequence at an infinite
hypernatural index~$n=H$.  Thus,
\begin{equation}
\label{21}
\lim_{n\to\infty} u_n = \text{st}(u_H).
\end{equation}
Here the standard part function ``st" associates to each finite
hyperreal, the unique finite real infinitely close to it (i.e., the
difference between them is infinitesimal).  This formalizes the
natural intuition that for ``very large" values of the index, the
terms in the sequence are ``very close" to the limit value of the
sequence.  Conversely, the standard part of a hyperreal~$u=[u_n]$
represented in the ultrapower construction by a Cauchy sequence
$(u_n)$, is simply the limit of that sequence:
\begin{equation}
\label{22}
\text{st}(u)=\lim_{n\to\infty}u_n.
\end{equation}
Formulas~\eqref{21} and \eqref{22} express limit and standard part in
terms of each other.  In this sense, the procedures of taking the
limit and taking the standard part are logically equivalent.






\subsection{Variable quantity}
\label{vc}
The mathematical term \megethos{} in ancient Greek has been translated
into Latin as \emph{quantitas}.  In modern languages it has two
competing counterparts: in English -- \emph{quantity},
\emph{magnitude};%
\footnote{\label{f23}The term ``magnitude'' is etymologically related
to \megethos.  Thus, \megethos{} in Greek and \emph{magnitudo} in
Latin both mean ``bigness''; ``big'' being \emph{mega}
($\mu\acute\varepsilon\gamma\alpha$) in Greek and \emph{magnum} in
Latin.}
in French -- \emph{quantit\'e}, \emph{grandeur}; in German --
\emph{Quantit\"{a}t}, \emph{Gr\"{o}sse}.  The term \emph{grandeur}
with the meaning \emph{real number} is still in use in (Bourbaki 1947
\cite{Bou}).  \emph{Variable quantity} was a primitive notion in
analysis as presented by Leibniz, l'H\^opital, and later Carnot and
Cauchy.  Other key notions of analysis were defined in terms of
variable quantities.  Thus, in Cauchy's terminology, a variable
quantity \emph{becomes} an \emph{infinitesimal} if it eventually drops
below any \emph{given} (i.e., constant) quantity (see Borovik \& Katz
\cite{BK} for a fuller discussion).  Cauchy notes that the
\emph{limit} of such a quantity is zero.  The notion of \emph{limit}
itself is defined as follows:
\begin{quote}
Lorsque les valeurs successivement attribu\'ees \`a une m\^eme
variable s'approchent ind\'efiniment d'une valeur fixe, de mani\`ere
\`a finir par en diff\'erer aussi peu que l'on voudra, cette
derni\`ere est appel\'ee la limite de toutes les autres (Cauchy,
\emph{Cours d'Analyse} \cite{Ca21}).
\end{quote}
Thus, Cauchy defined both infinitesimals and limits in
terms of the primitive notion of a variable quantity.  In Cauchy, any
variable quantity~$q$ that does not tend to infinity is expected to
decompose as the sum of a given quantity~$c$ and an
infinitesimal~$\alpha$:
\begin{equation}
\label{23}
q=c+\alpha.
\end{equation}
In his 1821 text \cite{Ca21}, Cauchy worked with a hierarchy of
infinitesimals defined by polynomials in a base infinitesimal
$\alpha$.  Each such infinitesimal decomposes as
\begin{equation}
\label{24}
\alpha^n(c+\varepsilon)
\end{equation}
for a suitable integer~$n$ and infinitesimal~$\varepsilon$.  Cauchy's
expression~\eqref{24} can be viewed as a generalisation of~\eqref{23}.

In Leibniz's terminology,~$c$ is an \emph{assignable} quantity
while~$\alpha$ and~$\varepsilon$ are \emph{inassignable}.  Leibniz's
\emph{transcendental law of homogeneity} (see entry~\ref{tlh})
authorized the replacement of the inassignable~$q=c+\alpha$ by the
assignable~$c$ since~$\alpha$ is negligible compared to~$c$:
\begin{equation}
q\adequal c
\end{equation}
(see entry~\ref{re} on relation~``$\adequal$'').  Leibniz emphasized
that he worked with a generalized notion of equality where expressions
were declared ``equal'' if they differed by a negligible term.
Leibniz's procedure was formalized in Robinson's B-approach by the
\emph{standard part function} (see entry~\ref{spf}), which assigns to
each finite hyperreal number, the unique real number to which it is
infinitely close.  As such, the standard part allows one to work
``internally'' (not in the technical NSA sense but) in the sense of
exploiting concepts already available in the toolkit of the historical
infinitesimal calculus, such as Fermat's \emph{adequality}
(entry~\ref{ad}), Leibniz's \emph{transcendental law of homogeneity}
(entry~\ref{tlh}), and Euler's \emph{principle of cancellation} (see
Bair et al.~\cite{B11}).  Meanwhile, in the A-approach as formalized
by Weierstrass, one is forced to work with ``external'' concepts such
as the multiple-quantifier~$\epsilon, \delta$ \emph{definitions} (see
entry~\ref{co}) which have no counterpart in the historical
infinitesimal calculus of Leibniz and Cauchy.

Thus, the notions of \emph{standard part} and \emph{epsilontic limit},
while logically equivalent (see entry~\ref{spf}), have the following
difference between them: the standard part principle corresponds to an
``internal'' development of the historical infinitesimal calculus,
whereas the epsilontic limit is ``external'' to it.


\subsection{Zeno's paradox of extension}
\label{ze}

Zeno of Elea (who lived about 2500 years ago) raised a puzzle (the
\emph{paradox of extension}, which is distinct from his better known
paradoxes of motion) in connection with treating any continuous
magnitude as though it consists of infinitely many indivisibles; see
(Sherry 1988 \cite{Sh88}); (Kirk et al.~1983 \cite{KRS}).  If the
indivisibles have no magnitude, then an \emph{extension} (such as
space or time) composed of them has no magnitude; but if the
indivisibles have some (finite) magnitude, then an extension composed
of them will be infinite.  There is a further puzzle: If a magnitude
is composed of indivisibles, then we ought to be able to add or
concatenate them in order to produce or increase a magnitude.  But
indivisibles are not next to one another; as limits or boundaries, any
pair of indivisibles is separated by what they limit.  Thus, the
concepts of addition or concatenation seem not to apply to
indivisibles.

The paradox need not apply to infinitesimals in Leibniz's sense,
however (see entry~\ref{36} on \emph{indivisibles and
infinitesimals}).  For, having neither zero nor finite magnitude,
infinitely many of them may be just what is needed to produce a finite
magnitude.  And in any case, the addition or concatenation of
infinitesimals (of the same dimension) is no more difficult to
conceive of than adding or concatenating finite magnitudes.  This is
especially important, because it allows one to apply arithmetic
operations to infinitesimals (see entry~\ref{lc} on the \emph{law of
continuity}).  See also (Reeder 2012 \cite{Re}).

\section*{Acknowledgments}

The work of Vladimir Kanovei was partially supported by RFBR grant
13-01-00006.  M.~Katz was partially funded by the Israel Science
Foundation grant no.~1517/12.  We are grateful to Reuben Hersh and
Martin Davis for helpful discussions, and to the anonymous referees
for a number of helpful suggestions.  The influence of Hilton Kramer
(1928-2012) is obvious.



\begin{thebibliography}{AM}




\bibitem{An11} Andersen, K.: One of Berkeley's arguments on
compensating errors in the calculus.  \emph{Historia Mathematica}
\textbf{38} (2011), no.~2, 219--231.


\bibitem{Arc} Archimedes, \emph{De Sphaera et Cylindro}, in
\emph{Archimedis Opera Omnia cum Commentariis Eutocii}, vol. I,
ed. J.L. Heiberg, B. G. Teubner, Leipzig, 1880.



\bibitem{Art} Arthur, R. (2007): Leibniz's syncategorematic
infinitesimals, Smooth Infinitesimal Analysis, and Newton's
Proposition 6.  See

http://www.humanities.mcmaster.ca/$\sim$rarthur/papers/LsiSiaNp6.rev.pdf


\bibitem{Bachet} Bachet: Diophanti Alexandrini, Arithmeticorum Liber~V
(Bachet's Latin translation).

\bibitem{B11} Bair, J.; B\l{}aszczyk, P.; Ely, R.; Henry, V.; Kanovei,
V.; Katz, K.; Katz, M.; Kutateladze, S.; McGaffey, T.; Schaps, D.;
Sherry, D.; Shnider, S.: Interpreting Euler's infinitesimal
mathematics (in preparation).

\bibitem{BAHE1} Bair J.; Henry V.: From Newton's fluxions to Virtual
Microscopes.  \emph{Teaching Mathematics and Computer Science}
\emph{V} (2007), 377-384.

\bibitem{BAHE2} Bair J.; Henry V.: From Mixed Angles to
Infinitesimals. \emph{The College Mathematics Journal} \emph{39}
(2008), no.~3, 230-233.

\bibitem{BH08} Bair J.; Henry V.: Analyse infinit\'esimale - Le
calculus red\'ecouvert.  Editions Academia Bruylant Louvain-la-Neuve
(Belgium) 2008, 189 pages 

D/2008/4910/33



\bibitem{Bel09} Bell, J. L.: Continuity and infinitesimals.  Stanford
Encyclopedia of philosophy.  Revised 20 july 2009.


\bibitem{Bi65A} Bishop, E.: Differentiable manifolds in complex
Euclidean space.  \emph{Duke Mathematical Journal} \textbf{32} (1965),
1--21.

\bibitem{Bi68} Bishop, E.: Mathematics as a numerical language.  1970
Intuitionism and Proof Theory (Proc. Conf., Buffalo, N.Y., 1968)
pp. 53--71.  North-Holland, Amsterdam.


\bibitem{Bi75} Bishop, E.: The crisis in contemporary mathematics.
Proceedings of the American Academy Workshop on the Evolution of
Modern Mathematics (Boston, Mass., 1974).  \emph{Historia Mathematica}
\textbf{2} (1975), no.~4, 507-517.


\bibitem{Bi77} Bishop, E.: Review: H. Jerome Keisler, Elementary
calculus, \emph{Bull. Amer. Math. Soc.} \textbf{83} (1977), 205-208.


\bibitem{Bi85} Bishop, E.: Schizophrenia in contemporary mathematics.
In Errett Bishop: reflections on him and his research (San Diego,
Calif., 1983), 1--32, \emph{Contemp. Math.} \textbf{39},
Amer. Math. Soc., Providence, RI, 1985 [originally distributed in
1973]

\bibitem{BKS} B\l aszczyk, P.; Katz, M.; Sherry, D.: Ten
misconceptions from the history of analysis and their debunking.
\emph{Foundations of Science} \textbf{18} (2013), no.~1, 43-74.  See
http://dx.doi.org/10.1007/s10699-012-9285-8

and http://arxiv.org/abs/1202.4153
%
%

\bibitem{BM} B\l{}aszczyk, P.; Mr\'owka, K.: Euklides, Elementy,
Ksi%
{e}gi V-VI. T\l{}umaczenie i komentarz [Euclid, Elements, Books V-VI.
Translation and commentary].  Copernicus Center Press, Krak\'ow, 2013.




\bibitem{BK} Borovik, A.; Katz, M.: Who gave you the
Cauchy--Weierstrass tale?  The dual history of rigorous calculus.
\emph{Foundations of Science}, \textbf{17} (2012), no.~3, 245--276.
See http://dx.doi.org/10.1007/s10699-011-9235-x

and http://arxiv.org/abs/1108.2885


\bibitem{Bos} Bos, H. J. M.: Differentials, higher-order differentials
and the derivative in the Leibnizian calculus.  \emph{Archive for
History of Exact Sciences} \textbf{14} (1974), 1--90.

\bibitem{Bos80} Bos, H.; Bunn, R.; Dauben, J.; Grattan-Guinness, I.;
Hawkins, T.; Pedersen, K.; From the calculus to set theory, 1630-1910.
An introductory history.  Edited by I. Grattan-Guinness.  Gerald
Duckworth \& Co. Ltd., London, 1980.


\bibitem{Bou} Bourbaki, N.: Th\'eorie de la mesure et de
l'int\'egration.  Introduction, Universit\'{e} Henri Poincar\'{e},
Nancy 1947.

\bibitem{Br} Br\aa{}ting, K.: A new look at E. G. Bj\"orling and the
Cauchy sum theorem.  \emph{Arch. Hist. Exact Sci.} \textbf{61} (2007),
no.~5, 519--535.
%
%



\bibitem{Bre94} Breger, H.: The mysteries of adaequare: a vindication
of Fermat.  \emph{Archive for History of Exact Sciences} \textbf{46}
(1994), no.~3, 193--219.


\bibitem{Ca21} Cauchy, A. L.: {\em Cours d'Analyse de L'Ecole Royale
Polytechnique.  Premi\`ere Partie.  Analyse alg\'ebrique\/}.  Paris:
Imprim\'erie Royale, 1821.  Online at

{\tiny http://books.google.com/books?id=\_mYVAAAAQAAJ\&dq=cauchy\&lr=\&source=gbs\_navlinks\_s}


\bibitem{Ca27} Cauchy, A.-L. (1815) Th\'eorie de la propagation des
ondes \`a la surface d'un fluide pesant d'une profondeur ind\'efinie
(published 1827, with additional Notes). Oeuvres, Series 1, Vol.~1,
4-318.


\bibitem{Ce07} Cerroni, C.: The contributions of Hilbert and Dehn to
non-Archimedean geometries and their impact on the Italian school.
\emph{ Revue d'Histoire des Math\'ematiques} \textbf{13} (2007),
no.~2, 259--299.

\bibitem{Ci} Cifoletti, G.: La m\'ethode de Fermat: son statut et sa
diffusion.  Alg\`ebre et comparaison de figures dans l'histoire de la
m\'ethode de Fermat.  \emph{Cahiers d'Histoire et de Philosophie des
Sciences. Nouvelle S\'erie} \textbf{33}. Soci\'et\'e Fran\c caise
d'Histoire des Sciences et des Techniques, Paris, 1990.


\bibitem{Co94} Connes, A.: Noncommutative geometry.  Academic Press,
Inc., San Diego, CA, 1994.


\bibitem{Da80} Dauben, J.: The development of the Cantorian set
theory.  In (Bos et al. 1980 \cite{Bos80}), pp.~181-219.

\bibitem{Da03} Dauben, J.: ``Abraham Robinson. 1918-1974,"
\emph{Biographical Memoirs of the National Academy of Sciences}
\textbf{82} (2003), 243-284.  Available at the address
http://www.nap.edu/html/biomems/arobinson.pdf and

http://www.nap.edu/catalog/10683.html

\bibitem{de} De Morgan, A.: On the early history of infinitesimals in
England.  \emph{Philosophical Magazine, Ser.~4} \textbf{4} (1852),
no.~26, 321--330.  See

http://www.tandfonline.com/doi/abs/10.1080/14786445208647134


\bibitem{Dir} Dirac, P.: The Principles of Quantum Mechanics.  4th
edition.  Oxford, at the Clarendon Press, 1958.



\bibitem{DM} Dossena, R.; Magnani, L.: Mathematics through diagrams:
microscopes in non-standard and smooth analysis.  \emph{Studies in
Computational Intelligence (SCI)} \textbf{64} (2007), 193--213.


\bibitem{Eh06} Ehrlich, P.: The rise of non-Archimedean mathematics
and the roots of a misconception. I. The emergence of non-Archimedean
systems of magnitudes.  \emph{Archive for History of Exact Sciences}
\textbf{60} (2006), no.~1, 1--121.


\bibitem{El12} Ely, R.: Loss of dimension in the history of calculus
and in student reasoning.  \emph{The Mathematics Enthusiast}
\textbf{9} (2012), no.~3, 303-326.

\bibitem{Euc} Euclid: Euclid's Elements of Geometry, edited, and
provided with a modern English translation, by Richard Fitzpatrick,
2007.  See

http://farside.ph.utexas.edu/euclid.html






\bibitem{ai} Euler, L.: Introductio in Analysin Infinitorum, Tomus
primus.  SPb and Lausana, 1748.

\bibitem{aie} Euler, L.: Introduction to analysis of the infinite.
Book I.  Translated from the Latin and with an introduction by John
D.~Blanton.  Springer-Verlag, New York, 1988 [translation of (Euler
1748 \cite{ai})]


\bibitem{ed} Euler, L.: Institutiones Calculi Differentialis.  SPb,
  1755.



\bibitem{Eu55} Euler, L.: Foundations of Differential Calculus.
English translation of Chapters 1--9 of (Euler 1755 \cite{ed}) by
D.\,Blanton, Springer, N.Y., 2000.


\bibitem{Fe04} Ferraro, G.: Differentials and differential
coefficients in the Eulerian foundations of the calculus.
\emph{Historia Mathematica} \textbf{31} (2004), no.~1, 34--61.


\bibitem{Fr71a} Freudenthal, H.: Cauchy, Augustin-Louis.  In
Dictionary of Scientific Biography, ed. by C. C. Gillispie, vol. 3
(New York: Charles Scribner's sons, 1971), 131-148.


\bibitem{Ge46} Gerhardt:, C. I. (ed.): Historia et Origo calculi
differentialis a G. G. Leibnitio conscripta, ed. C. I. Gerhardt,
Hannover, 1846.

\bibitem{Ge50} Gerhardt, C. I. (ed.): Leibnizens mathematische
Schriften (Berlin and Halle: Eidmann, 1850-1863).
%
%


\bibitem{Gran} Granger, G.: Philosophie et math\'ematique
leibniziennes.  \emph{Revue de M\'eta\-physique et de Morale 86e
Ann\'ee}, No. 1 (Janvier-Mars 1981), pp.~1-37.

\bibitem{Gras} Grassmann, Hermann: \emph{Lehrbuch der Arithmetik},
Enslin, Berlin 1861.



\bibitem{Gr08a} Gray, J.: Plato's ghost. The modernist transformation
of mathematics.  Princeton University Press, Princeton, NJ, 2008.

\bibitem{Gr08b} Gray, J.: A short life of Euler.  \emph{BSHM
Bulletin. Journal of the British Society for the History of
Mathematics} \textbf{23} (2008), no.~1, 1--12.

\bibitem{Hal} Hall, Rupert; Hall, Marie Boas, eds.  Unpublished
Scientific Papers of Isaac Newton, Cambridge University Press, 1962,
p.15-19, 31-37, quoted in I. Bernard Cohen, Richard S. Westfall,
Newton, Norton Critical Edition, 1995, p. 377-386.

\bibitem{Ha07} Hawking, S., Ed.: The mathematical breakthroughs that
changed history. Running Press, Philadelphia, PA, 2007.

\bibitem{Han} Hankel, Hermann: Zur Geschichte der Mathematik in
Alterthum und Mittelalter.  Teubner, Leipzig 1876.

\bibitem{He97} Heath, T. (ed.) The works of Archimedes.  Cambridge
University Press, Cambridge, 1897.

\bibitem{Hew48} Hewitt, E.: Rings of real-valued continuous
functions. I.  \emph{Trans. Amer. Math. Soc.} \textbf{64} (1948),
45--99.

\bibitem{Hei} Heiberg, J. (ed.): Archimedis Opera Omnia cum
Commentariis Eutocii, vol.~I.  Teubner, Leipzig 1880.


\bibitem{He73} Heijting, A.: Address to Professor A. Robinson.  At the
occasion of the Brouwer memorial lecture given by Prof. A.~Robinson on
the 26th April 1973.  \emph{Nieuw Arch. Wisk. (3)} \textbf{21} (1973),
134--137.  MathScinet Review at

http://www.ams.org/mathscinet-getitem?mr=434756


\bibitem{Her} Herzberg, F.: Stochastic calculus with infinitesimals.
Lecture Notes in Mathematics, 2067. Springer, Heidelberg, 2013.


\bibitem{Hi} Hilbert, D.: Grundlagen der Geometrie.  Festschrift zur
Enth\"{u}llung des Gauss-Weber Denkmals in G\"{o}ttingen, Leipzig
1899.


\bibitem{Hi00} Hilbert, D.: Les principes fondamentaux de la
g\'eom\'etrie.  \emph{Annales scientifiques de l'E.N.S. ~3$^e$
s\'erie} \textbf{17} (1900), 103-209.



\bibitem{Ho} H\"{o}lder, O.: Die Axiome der Quantit\"{a}t und die
Lehre vom Mass.  Berichte \"{u}ber die Verhandlungen der K\"{o}niglich
S\"{a}chsischen Gesellschaft der Wissenschaften zu Leipzig.
Mathematisch-Physische Classe, 53, Leipzig 1901, 1--63.


\bibitem{Is} Ishiguro, H.: Leibniz's philosophy of logic and language.
Second edition.  Cambridge University Press, Cambridge, 1990.


\bibitem{Je11} Jesseph, D.: Leibniz on the Elimination of
infinitesimals: Strategies for finding truth in fiction, 27 pages.  In
Leibniz on the interrelations between mathematics and philosophy,
edited by Norma B. Goethe, Philip Beeley and David Rabouin, Archimedes
Series, Springer Verlag, 2013.

\bibitem{kan} Kanovei, V.: The correctness of Euler's method for the
factorization of the sine function into an infinite product,
\emph{Russian Mathematical Surveys} \textbf{43} (1988), 65--94.


\bibitem{KKM} Kanovei, V.; Katz, M.; Mormann, T.: Tools, Objects, and
Chimeras: Connes on the Role of Hyperreals in Mathematics.
\emph{Foundations of Science} \textbf{18} (2013), no.~2, 259--296.
See http://dx.doi.org/10.1007/s10699-012-9316-5

and http://arxiv.org/abs/1211.0244


\bibitem{kr} Kanovei, V.; Reeken, M.: Nonstandard analysis,
axiomatically.  Springer Monographs in Mathematics, Berlin: Springer,
2004.


\bibitem{KS} Kanovei, V.; Shelah, S.: A definable nonstandard model of
the reals.  \emph{Journal of Symbolic Logic} \textbf{69} (2004),
no.~1, 159--164.



\bibitem{KK1} Katz, K.; Katz, M.: Zooming in on infinitesimal~$1-.9..$
in a post-triumvirate era.  \emph{Educational Studies in Mathematics}
\textbf{74} (2010), no.~3, 259-273.  

See http://arxiv.org/abs/arXiv:1003.1501


\bibitem{KK2} Katz, K.; Katz, M.: When is $.999\ldots$ less than $1$?
\emph{The Montana Mathematics Enthusiast} \textbf{7} (2010), No.~1,
3--30.


\bibitem{KK11b} Katz, K.; Katz, M.: Cauchy's continuum.
\emph{Perspectives on Science} \textbf{19} (2011), no.~4, 426-452.
See http://www.mitpressjournals.org/toc/posc/19/4

and http://arxiv.org/abs/1108.4201


\bibitem{KK11d} Katz, K.; Katz, M.: Meaning in classical mathematics:
is it at odds with Intuitionism?  \emph{Intellectica}, \textbf{56}
(2011), no.~2, 223-302, see 

http://arxiv.org/abs/1110.5456


\bibitem{KK11a} Katz, K.; Katz, M.: A Burgessian critique of
nominalistic tendencies in contemporary mathematics and its
historiography.  \emph{Foundations of Science} \textbf{17} (2012),
no.~1, 51--89.  See http://dx.doi.org/10.1007/s10699-011-9223-1

and http://arxiv.org/abs/1104.0375

%

\bibitem{KL} Katz, M.; Leichtnam, E.: Commuting and noncommuting
infinitesimals.  \emph{American Mathematical Monthly} \textbf{120}
(2013), no. 7, 631--641.

See http://arxiv.org/abs/1304.0583

\bibitem{KSS13} Katz, M.; Schaps, D.; Shnider, S.: Almost Equal: The
Method of Adequality from Diophantus to Fermat and Beyond.
\emph{Perspectives on Science} \textbf{21{}} (2013), no.~3, 283-324.
See http://arxiv.org/abs/1210.7750


\bibitem{KS2} Katz, M.; Sherry, D.: Leibniz's laws of continuity and
homogeneity.  \emph{Notices of the American Mathematical Society}
\textbf{59} (2012), no.~11, 1550-1558.  See
http://www.ams.org/notices/201211/ and http://arxiv.org/abs/1211.7188


\bibitem{KS1} Katz, M.; Sherry, D.: Leibniz's infinitesimals: Their
fictionality, their modern implementations, and their foes from
Berkeley to Russell and beyond.  \emph{Erkenntnis} \textbf{78} (2013),
no.~3, 571--625.  See http://dx.doi.org/10.1007/s10670-012-9370-y and
http://arxiv.org/abs/1205.0174






\bibitem{KT2} Katz, M.; Tall, D.: A Cauchy-Dirac delta function.
\emph{Foundations of Science} \textbf{18} (2013), no.~1, 107-123.  See
http://dx.doi.org/10.1007/s10699-012-9289-4 and
http://arxiv.org/abs/1206.0119



\bibitem{Ke08} Keisler, H. J.: The ultraproduct
construction. Ultrafilters across mathematics, 163-179,
\emph{Contemp. Math.}, \textbf{530}, Amer. Math. Soc., Providence, RI,
2010.

\bibitem{KRS} Kirk, G.S.; Raven, J.E.; Schofield M.: Presocratic
Philosophers.  Cambridge UP, Cambridge 1983, \S 316.


\bibitem{Kl08} Klein, F.: Elementary Mathematics from an Advanced
Standpoint.  Vol. I.  Arithmetic, Algebra, Analysis.  Translation by
E. R. Hedrick and C. A. Noble [Macmillan, New York, 1932] from the
third German edition [Springer, Berlin, 1924].  Originally published
as Elementarmathematik vom h\"oheren Standpunkte aus (Leipzig, 1908).


\bibitem{Kl80} Kline, M.: Mathematics. The loss of certainty. Oxford
University Press, New York, 1980.





\bibitem{Lau87} Laugwitz, D.: Hidden lemmas in the early history of
infinite series.  \emph{Aequationes Mathematicae} \textbf{34} (1987),
264--276.


\bibitem{Lau89} Laugwitz, D.: Definite values of infinite sums:
aspects of the foundations of infinitesimal analysis around 1820.
\emph{Archive for History of Exact Sciences} \textbf{39} (1989),
no.~3, 195--245.


\bibitem{Lau92} Laugwitz, D.: Early delta functions and the use of
infinitesimals in research.  \emph{Revue d'histoire des sciences}
\textbf{45} (1992), no.~1, 115--128.



\bibitem{Le01c} Leibniz, G. (1701) \emph{Cum Prodiisset}\ldots mss
``Cum prodiisset atque increbuisset Analysis mea infinitesimalis ..."
in Gerhardt \cite{Ge46}, pp.~39--50.


\bibitem{Le02} Leibniz, G. (1702) Letter to Varignon, 2 feb.~1702, in
Gerhardt \cite[vol.~IV, pp. 91--95]{Ge50}.


\bibitem{Le10b} Leibniz, G. (1710) Symbolismus memorabilis calculi
algebraici et infinitesimalis in comparatione potentiarum et
differentiarum, et de lege homogeneorum transcendentali.  In Gerhardt
\cite[vol.~V, pp.~377-382]{Ge50}.

\bibitem{Le89} Leibniz, G.: La naissance du calcul diff\'erentiel.  26
articles des Acta Eruditorum.  Translated from the Latin and with an
introduction and notes by Marc Parmentier.  With a preface by Michel
Serres.  Mathesis.  Librairie Philosophique J. Vrin, Paris, 1989.

\bibitem{Le1993} Leibniz, G.: De quadratura arithmetica circuli
ellipseos et hyperbolae cujus corollarium est trigonometria sine
tabulis.  Edited, annotated and with a foreword in German by Eberhard
Knobloch.  Abhandlungen der Akademie der Wissenschaften in
G\"ottingen.  Mathematisch-Physikalische Klasse. Folge 3 [Papers of
the Academy of Sciences in G\"ottingen.  Mathematical-Physical
Class. Series 3], 43. Vandenhoeck \& Ruprecht, G\"ottingen, 1993.





\bibitem{Lo} {\L}o{\'s}, J.: Quelques remarques, th\'eor\`emes et
probl\`emes sur les classes d\'efi\-nissables d'alg\`ebres.  In
Mathematical interpretation of formal systems, {98--113},
North-Holland Publishing Co., {Amsterdam}, {1955}.
%
%


\bibitem{lux} Luxemburg, W.: What is nonstandard analysis?  Papers in
the foundations of mathematics.  \emph{American Mathematical Monthly}
\textbf{80} (1973), no.~6, part II, 38--67.


\bibitem{MD} Magnani, L.; Dossena, R.: Perceiving the infinite and the
infinitesimal world: unveiling and optical diagrams in mathematics.
\emph{Foundations of Science} \textbf{10} (2005), no.~1, 7--23.


\bibitem{Ma96} Mancosu, P.: Philosophy of mathematics and mathematical
practice in the seventeenth century.  The Clarendon Press, Oxford
University Press, New York, 1996.







\bibitem{Mc23} McClenon, R.: A Contribution of Leibniz to the History
of Complex Numbers.  \emph{American Mathematical Monthly} \textbf{30}
(1923), no.~7, 369-374.

\bibitem{mt} McKinzie, M.; Tuckey, C.: Hidden lemmas in Euler's
summation of the reciprocals of the squares, \emph{Archive for History
of Exact Sciences} \textbf{51} (1997), 29--57.

\bibitem{Me} Meschkowski, H.: Aus den Briefbuchern Georg Cantors.
\emph{Archive for History of Exact Sciences} \textbf{2} (1965),
503-519.


\bibitem{MK} Mormann, T.; Katz, M.: Infinitesimals as an issue of
neo-Kantian philosophy of science.  \emph{HOPOS: The Journal of the
International Society for the History of Philosophy of Science}, to
appear.



\bibitem{Mu11} Mumford, D.: Intuition and rigor and Enriques's
quest. \emph{Notices Amer. Math. Soc.} \textbf{58} (2011), no.~2,
250--260.


\bibitem{Ne} Newton, I. (1671), \emph{Methodus Fluxionum}.  English
version \emph{The method of fluxions and infinite series} (1736).

\bibitem{NEM} N\'u\~nez, R.; Edwards, L.; Matos, J.: Embodied
cognition as grounding for situatedness and context in mathematics
education.  {\em Educational Studies in Mathematics\/} 39 (1999),
no.~1-3, 45-65, DOI: 10.1023/A:1003759711966.
%
%


\bibitem{PK} Parkhurst, W.; Kingsland, W.: Infinity and the
infinitesimal.  \emph{The Monist} \textbf{35} (1925), 633-666.


\bibitem{Pr} Proietti, C.: Natural Numbers and Infinitesimals: A
Discussion between Benno Kerry and Georg Cantor.  \emph{History and
Philosophy of Logic} \textbf{29} (2008), no.~4, 343--359.


\bibitem{Re} Reeder, P.: Infinitesimals for Metaphysics: Consequences
for the Ontologies of Space and Time.  Degree Doctor of Philosophy,
Ohio State University, Philosophy, 2012.


\bibitem{Ro61} Robinson, A.: Non-standard analysis.
\emph{Nederl. Akad. Wetensch. Proc. Ser. A} \textbf{64} =
\emph{Indag. Math.} \textbf{23} (1961), 432--440 [reprinted in
Selected Works, see (Robinson 1979~\cite[pp.~3--11]{Ro79}).



\bibitem{Ro68} Robinson, A.: ``Reviews: Foundations of Constructive
Analysis''.  \emph{American Mathematical Monthly} \textbf{75} (1968),
no.~8, 920--921.

\bibitem{Ro79} Robinson, A.: Selected papers of Abraham Robinson.
Vol. II.  Nonstandard analysis and philosophy.  Edited and with
introductions by W. A. J. Luxemburg and S. K\"orner.  Yale University
Press, New Haven, Conn., 1979.

\bibitem{Ru03} Russell, B.: The principles of mathematics.  Vol. I.
Cambridge Univ. Press, Cambridge, 1903.

\bibitem{She87} Sherry, D.: The wake of Berkeley's Analyst:
\emph{rigor mathematicae}? \emph{Stud. Hist. Philos. Sci.}
\textbf{18} (1987), no.~4, 455--480.

\bibitem{Sh88} Sherry, D.: Zeno's metrical paradox revisited.
\emph{Philosophy of Science} \textbf{55} (1988), no.~1, 58--73.

\bibitem{SK} Sherry, D.; Katz, M.: Infinitesimals, imaginaries,
ideals, and fictions.  \emph{Studia Leibnitiana} (2013), to appear.
See http://arxiv.org/abs/1304.2137

\bibitem{Si} Simson, R.: The Elements of Euclid viz. the first six
Books, together with the eleventh and twelfth.  The Errors by which
Theon and others have long ago vitiated these Books are corrected and
some of Euclid's Demonstrations are restored.  Robert \& Andrew
Foulis, Glasgow, 1762.


\bibitem{Sk33} Skolem, T.: \"Uber die Unm\"oglichkeit einer
vollst\"andigen Charakterisierung der Zahlenreihe mittels eines
endlichen Axiomensystems.  \emph{Norsk Mat. Forenings Skr., II. Ser.}
No.~1/12 (1933), 73-82.

\bibitem{Sk34} Skolem, T.: \"Uber die Nicht-charakterisierbarkeit der
Zahlenreihe mittels endlich oder abz\"ahlbar unendlich vieler Aussagen
mit ausschliesslich Zahlenvariablen.  \emph{Fundamenta Mathematicae}
\textbf{23} (1934), 150-161.

\bibitem{Sk55} Skolem, T.: Peano's axioms and models of arithmetic.
In Mathematical interpretation of formal systems, pp. 1--14.
North-Holland Publishing Co., Amsterdam, 1955.

\bibitem{So} Solovay, R.: A model of set-theory in which every set of
reals is Lebesgue measurable.  \emph{Annals of Mathematics (2)}
\textbf{92} (1970), 1--56.


\bibitem{Sto} Stolz, O.: \emph{Vorlesungen \"{u}ber Allgemeine
Arithmetik}, Teubner, Leipzig 1885.


\bibitem{TK} Tall, D.; Katz, M.: A cognitive analysis of Cauchy's
conceptions of function, continuity, limit, and infinitesimal, with
implications for teaching the calculus.  \emph{Educational Studies in
Mathematics} (2013), to appear.

\bibitem{Th} Tho, T.: Equivocation in the foundations of Leibniz's
infinitesimal fictions. \emph{Society and Politics} \textbf{6} (2012),
No.~2 (12).


\bibitem{Vi13} Vickers, P.: Understanding Inconsistent Science.
Oxford University Press, Oxford, 2013.


\bibitem{Wa} Wallis, J. (1656) The arithmetic of infinitesimals.
Translated from the Latin and with an introduction by Jaequeline
A. Stedall.  \emph{Sources and Studies in the History of Mathematics
and Physical Sciences}.  Springer-Verlag, New York, 2004.







\end{thebibliography}
\end{document}